\documentclass[12pt,reqno]{amsart}

\usepackage[LGR,T1]{fontenc}
\usepackage[utf8]{inputenc}

\title[Time reversal of diffusion processes]{Time reversal of diffusion processes under a finite entropy condition}
\date{August, 2022}

\author[Cattiaux]{Patrick Cattiaux}
\address{Patrick Cattiaux. Institut de Mathématiques de Toulouse.  Université Paul Sabatier, Toulouse,   France.}
\email{patrick.cattiaux@math.univ-toulouse.fr}

\author[Conforti]{Giovanni Conforti}
\address{Giovanni Conforti, Département de Mathématiques Appliquées, \'Ecole Polytechnique, Palaiseau, France.}
\email{giovanni.conforti@polytechnique.edu}

\author[Gentil]{Ivan Gentil}
\address{Ivan Gentil. Institut Camille Jordan, Université Claude Bernard, Lyon, France. }
\email{gentil@math.univ-lyon1.fr}

\author[Léonard]{Christian Léonard}
\address{Christian Léonard. Modal-X,  Université Paris Nanterre, France.}
\email{christian.leonard@math.cnrs.fr}

\thanks{This research is partially  granted by the projects  EFI (ANR-17-CE40-003), SPOT (ANR-20-CE40-0014) and  Labex MME-DII (ANR-11-LBX-0023)}

\keywords{Time-reversal, diffusion process, stochastic derivative, relative entropy, random walk,  entropic optimal transport}
\subjclass[2010]{60J60, 60J25}

\usepackage[english]{babel}

\usepackage{amssymb, amsmath, amsfonts, latexsym, enumerate}
\usepackage[usenames,dvipsnames]{pstricks}
\usepackage{epsfig}
\RequirePackage[colorlinks,linkcolor=blue,citecolor=blue,urlcolor=blue]{hyperref}

\newtheorem{theorem}[equation]{Theorem}
\newtheorem{lemma}[equation]{Lemma}
\newtheorem{proposition}[equation]{Proposition}
\newtheorem{corollary}[equation]{Corollary}
\newtheorem{claim}[equation]{Claim}

\newtheorem{definition}[equation]{Definition}
\newtheorem{definitions}[equation]{Definitions}

\newtheorem{hypotheses}[equation]{Hypotheses}
\newtheorem{hypothesis}[equation]{Hypothesis}

\theoremstyle{remark}
\newtheorem{remark}[equation]{Remark}
\newtheorem{remarks}[equation]{Remarks}

\numberwithin{equation}{section}

\newcommand{\RR}{\mathbb{R}}
\newcommand{\Rn}{\mathbb{R}^n}

\newcommand{\1}{\mathbf{1}}

\newcommand{\ttimes}{\!\times\!}

\newcommand\pf{_{\#}}

\newcommand{\Leb}{\mathrm{Leb}}

\renewcommand{\ae}{\textrm{-}\mathrm{a.e.}}
\newcommand{\Id}{\mathrm{Id}}
\newcommand{\scal}{\!\cdot\!}

\DeclareMathOperator{\dom}{dom}

\DeclareMathOperator{\supp}{supp}

\DeclareMathOperator{\dive}{div}
\DeclareMathOperator{\tr}{tr}
\DeclareMathOperator{\MP}{MP}

\newcommand{\boulette}[1]{$\bullet$\ Proof of #1.}
\newcommand{\Boulette}[1]{\par\medskip\noindent $\bullet$\ Proof of #1.}
\newcommand{\sbt}{\,\begin{picture}(-1,1)(-1,-3)\circle*{3}\end{picture}\ }

\newcommand\Lim[1]{\lim_{#1\rightarrow\infty}}

\newcommand\Limh{\lim_{h\to 0^+}}

\newcommand\II[2]{\int_{[#1,#2]}}

\newcommand{\cadlag}{c\`adl\`ag}
\newcommand{\ud}{\frac{1}{2}}

\newcommand{\ii}{[0,T]}
\newcommand\XX{ \mathcal{X}}
\newcommand\XXX{\XX^2}
\newcommand\XXb{{\overline{\XX}}}
\newcommand\iX{{\ii\ttimes\XX}}
\newcommand\iXX{{\ii\times\XXX}}
\newcommand{\ZZ}{\Rn}

\newcommand{\iZ}{\ii\times\ZZ}

\newcommand\OO{\Omega}

\newcommand{\UU}{ \mathcal{U}}
\newcommand{\UP}{ \UU^P}
\newcommand{\UQ}{ \UU^Q}
\newcommand{\UPb}{ \UU ^{ \Pb}}

\newcommand\PX{\mathrm{P}(\XX)}

\newcommand\MX{\mathrm{M}(\XX)}

\newcommand\PZ{\mathrm{P}(\ZZ)}

\newcommand\PO{\mathrm{P}(\OO)}
\newcommand\MO{\mathrm{M}(\OO)}
\newcommand{\CcZ}{C^2_c(\ZZ)}

\newcommand\Iii{\int_{\ii}}
\newcommand\IX{\int_{\XX}}
\newcommand\IXX{\int_{\XXX}}
\newcommand\IiX{\int_{\iX}}
\newcommand\IiXX{\int_{\ii\times\XXX}}
\newcommand\IZ{\int_{\ZZ}}

\newcommand\IiZ{\int_{\iZ}}

\newcommand\Xb{\overline{X}}
\newcommand\Pb{\overline{P}}

\newcommand\Qb{\overline{Q}}

\newcommand{\Lf}{\overrightarrow{L}}
\newcommand{\Lb}{\overleftarrow{L}}
\newcommand{\LL}{ \mathcal{L}}
\newcommand{\LLf}{\overrightarrow{ \mathcal{L}}}
\newcommand{\LLb}{\overleftarrow{ \mathcal{L}}}
\newcommand{\LLc}[1]{ \LL^{ \mathrm{cu},#1}}
\newcommand{\LLo}[1]{ \LL^{ \mathrm{os},#1}}

\newcommand{\Gaf}{\overrightarrow{\Gamma}}
\newcommand{\Gab}{\overleftarrow{\Gamma}}
\newcommand{\Jf}{\overrightarrow{J}}

\newcommand{\jf}{\overrightarrow{j}}
\newcommand{\jb}{\overleftarrow{j}}

\newcommand{\vf}{\overrightarrow{ \vv}}
\newcommand{\vb}{\overleftarrow{\vv}}
\newcommand{\vam}{\vv ^{ \aa,\mm}}
\newcommand{\vc}[1]{\vv ^{ \mathrm{cu},#1}}
\newcommand{\vo}[1]{\vv ^{ \mathrm{os},#1}}
\newcommand{\vcu}{\vv ^ \mathrm{cu}}
\newcommand{\vos}{\vv ^ \mathrm{os}}
\renewcommand{\bf}{\overrightarrow{\beta}}
\newcommand{\bb}{\overleftarrow{\beta}}
\newcommand{\bc}[1]{\beta ^{ \mathrm{cu},#1}}
\newcommand{\bo}[1]{\beta ^{ \mathrm{os},#1}}
\newcommand{\bcu}{\beta ^ \mathrm{cu}}
\newcommand{\bos}{\beta ^ \mathrm{os}}

\renewcommand{\aa}{ \mathsf{a}}
\renewcommand{\gg}{ \mathsf{g}}
\newcommand{\mm}{ \mathsf{m}}

\newcommand{\pp}{ \mathsf{p}}
\newcommand{\qq}{ \mathsf{q}}
\newcommand{\rr}{ \mathsf{r}}
\newcommand{\vv}{ \mathsf{v}}
\newcommand{\hh}{ \mathsf{h}}

\newcommand{\pb}{\bar{\pp}}
\newcommand{\qb}{\bar{\qq}}

\newcommand{\bbb}[2]{ \beta ^{ #1|#2}}

\newcommand{\sy}{\sum _{y:y\sim x}}
\newcommand{\sx}{\sum _{x\in\XX}}
\newcommand{\sxy}{\sum _{ (x,y):x\sim y}}

\begin{document}

\maketitle 

\begin{abstract} 
Motivated by entropic optimal transport, time reversal of diffusion processes is revisited. 
An  integration by parts formula is derived for the carré du champ of a  Markov process in an abstract space. It leads to a time reversal formula for a wide class of diffusion processes in $\ZZ$ possibly with singular drifts, extending the  already known results in this domain. 

The proof of  the integration by parts formula relies on stochastic derivatives. This formula is applied to compute the semimartingale characteristics of the time-reversed $P^*$ of a  diffusion measure $P$ provided that the relative entropy of $P$ with respect to another diffusion measure $R$ is finite, and the semimartingale characteristics of the time-reversed $R^*$ are known (for instance when the reference path measure $R$ is reversible). 
 
 As an illustration of the robustness of this method, the integration by parts formula is also employed to derive a time-reversal formula for a random walk on a graph.

\end{abstract}

\tableofcontents

\section{Introduction}

The time-reversed $(Y^*_t:=Y _{ T-t}, 0\le t\le T)$ of a Markov process $(Y_t, 0\le t\le T)$ remains a Markov process. Consequently, the problem of finding its Markov generator  arises naturally. The answer to this problem is given by the so-called time reversal formula.
More precisely, we  establish   an integration by parts  formula connecting the carré du champ of a Markov process (its Dirichlet form) with its backward and forward generators, see \eqref{eq-24-intro}. This result extends the well known case of reversible processes where forward and backward generators are equal, and is valid under mild regularity assumptions. 
Then, we apply this abstract integration by parts formula (IbP) to identify the semimartingale characteristics of a time-reversed diffusion process, see Theorems \ref{res-08-intro} and  \ref{res-08-extension}. 

It is worth mentioning that this IbP formula  allows a unified study of time reversal for diffusion processes and processes with jumps. See Section \ref{sec-rw} below  where time reversal of random walks on graphs is investigated to demonstrate the robustness of this strategy, and  the companion paper \cite{CL20} where  the time-reversal of a Markov process with jumps is investigated.

During the eighties, Föllmer   gave  a smart   proof of the time reversal formula for a diffusion process driven by a Brownian motion under  a finite entropy condition, using Nelson stochastic derivatives  \cite{Foe85b,Foe86}.
In the present article, we revisit F\"ollmer's proof working out in more detail some technical steps and  extending it to more general diffusion processes. In doing so, we keep its powerful guideline based on stochastic derivatives and entropic arguments. 

We stress that our version of the time reversal formula, as well as F\"ollmer's one, holds under  a  finite entropy  hypothesis implying a  low regularity of the drift field, $L^2$ being typically enough. Such a situation is not covered by the other main results in the field.

Besides being an interesting topic in its own right,  last years have seen a renewed interest in  time reversal  because of its applications to the Schr\"odinger problem (a.k.a.\ entropic optimal transport), see \eqref{eq-Sch-pb} below, and functional inequalities. To give some examples, in \cite{CGP14} and \cite{GLR17a} a fluid-dynamic (Benamou-Brenier) formulation of entropic optimal transport is derived leveraging time reversal arguments, and similar ideas are also used in \cite{BCGL19} in a mean field setting. In all these applications, it is of fundamental importance of having a result applicable to diffusions whose drift is only square integrable, as no more than this can be expected assuming only finite entropy with respect to the Brownian motion. We refer to  subsection \emph{``Entropic and deterministic optimal transports''} below for a slightly more accurate discussion of the links between time reversal and the Schr\"odinger problem.

Regarding functional inequalities, it is worth mentioning that Fontbona and Jourdain \cite{FJ16} recover and extend the Bakry-\'Emery criterion using an approach based on time reversal. Also using time reversal in a crucial manner, a simple proof of the logarithmic Sobolev inequality is proposed by Léonard in \cite{Leo12d},  Gentil, Léonard, Ripani and Tamanini \cite{GLR17b} derive the HWI inequality, and
 Karatzas, Schachermayer and Tschiderer \cite{KST18}   obtain   pathwise results about the exponential rate of convergence to equilibrium of some Wasserstein gradient flows and another   proof of the  HWI inequality.

 \subsection*{Outline of the article}  Next Section \ref{sec-sd} gathers basic notions about Nelson stochastic derivatives that will be used throughout the paper. Related technical results which are necessary during our proofs are postponed to the appendix Section \ref{sec-leo}. The main general result of the article is the integration by parts formula stated at Theorem \ref{res-02}. Section \ref{sec-df} is devoted to its proof. The time reversal formulas that  we obtain for diffusion processes    in Sections \ref{sec-diff1} and \ref{sec-diff2},  and random walks in Section \ref{sec-rw}   are corollaries of this theorem.   These   time reversal formulas are stated at Theorems  \ref{res-08}, \ref{res-08-extension}, \ref{res-20} and \ref{res-21}.  Finally, in  Section \ref{sec-cod},  the current-osmosis decomposition of an entropic interpolation in a diffusion setting is discussed in detail to illustrate our main motivation for revisiting time reversal under a finite entropy condition.

Theorem  \ref{res-08-extension} which is stated in this introductory section is an extension of Theorem  \ref{res-08}. Its proof is almost verbatim the same as Theorem  \ref{res-08}'s proof.

\subsection*{Entropic and deterministic optimal transports}

Let us start presenting some aspect of our main motivation for revisiting time reversal of Markov processes.
\\
Entropic optimal transport dates back to the seminal articles \cite{Sch31,Sch32} by Schrödinger and was rigorously rephrased in terms of large deviations of empirical measures of particle systems by Föllmer in his Saint-Flour lecture notes \cite{Foe85}.  One wants to minimize  the relative entropy 
\begin{align*}
H(P|R):=E_P\log(dP/dR)
\end{align*}
 with respect to the law $R$ of some reference Markov process on a time interval $[0,T]$ among all Markov  measures $P$ with prescribed initial and final marginals:
 \begin{align}\label{eq-Sch-pb}
 \inf \{ H(P|R);\ &P:  P_0= \mu_0, P_T= \mu_T\}.
 \end{align}
 
To fix the ideas in this introductory discussion,  following Schrödinger and Föllmer, our reference measure $R$ is the law of a Brownian motion. We denote by  $P_t$ the $t$-marginal of $P$, i.e.\ the law under $P$ of the position at time $t$, and $P^*$ the time-reversal of $P.$   The quantity $H(P|R)-H(P_0|R_0)$ appears as an average \emph{forward} kinetic action (again a result by Föllmer). Since time-reversal is a one-one mapping, we have 
\begin{align}\label{eq-38}
H(P^*|R^*)=H(P|R)
\end{align}
 which allows us to interpret $H(P|R)-H(P_T|R_T)$  as a \emph{backward} kinetic action. Taking the half sum, we arrive at 
 \begin{align}\label{eq-35a}
H(P|R)= \mathrm{function}(P_0,P_T)
 	+ A _{ \mathrm{cu}}(P)
		+A  _{ \mathrm{os}}( P),
 \end{align}
where the current action term $A _{ \mathrm{cu}}(P)$ is  purely kinetic  with a direct interpretation in terms of \emph{deterministic optimal transport}:
\begin{align*}
\inf \left\{A _{ \mathrm{cu}}(P);\ P:P_0= \mu_0,P_T= \mu_T)\right\} 
=T ^{ -1}W_2^2( \mu_0, \mu_T),
\end{align*}
with $W_2$ the standard quadratic Wasserstein distance. This is the Benamou-Brenier formula.  It turns out that the osmotic action term
\begin{align}\label{eq-35b}
A  _{ \mathrm{os}}( P)= \mathsf{A} _{ \mathrm{os}}( [P])
\end{align}
 only depends on the marginal flow $ [P]:=(P_t) _{ 0\le t\le T}$ of $P$ (it is directly linked to some Fisher information). This important identity follows from 
 the \emph{time reversal formula}, which is the main goal of this paper. 
  
  The decomposition \eqref{eq-35a} with \eqref{eq-35b} plays a major rôle in the comparison between  deterministic  and  entropic optimal transports.  In particular, we see that for a given flow of marginals $ \mu:=( \mu_t) _{ 0\le t\le T},$
\begin{align}\label{Lag-action}
\begin{split}
\mathsf{A}(\mu):=\inf \{ H(P|R);\ &P:  P_t= \mu_t, 0\le t\le T\} \\
&=\mathrm{function}( \mu_0, \mu_T)
+ \mathsf{A}_{ \mathrm{BB}}( \dot\mu)
 			+ \mathsf{A}  _{ \mathrm{os}}(  \mu)
\end{split}
\end{align}
where $ \mathsf{A} _{ \mathrm{BB}}(\dot \mu):=\inf \left\{A _{ \mathrm{cu}}(P);\ P:P_t= \mu_t, 0\le t\le T\right\}$ is the Benamou-Brenier action of  $ \mu:$   the fundamental notion of Otto calculus on the Wasserstein space of probability measures leading to the definition of the tangent vector $\dot\mu,$ see \cite{AGS05,Vill09}. \emph{In view of \eqref{Lag-action}, the osmotic action, whose appearance is tightly connected to time reversal, quantifies the difference between McCann displacement interpolations (attached to the standard deterministic quadratic transport)  and their entropic approximations, i.e.\ marginal flows of solutions of the Schrödinger problem \eqref{eq-Sch-pb} (attached to Brownian path measures $R$ with a diffusion coefficient tending to zero)}. 

In addition, at least formally, we see with \eqref{Lag-action} that the Schrödinger problem \eqref{eq-Sch-pb} is recast as a standard least action principle on the Wasserstein space of probability measures, where $ \mathsf{A}_{ \mathrm{BB}}(\dot\mu)$ is the kinetic action and $-  \mathsf{A}  _{ \mathrm{os}}(  \mu)$  is the action of some scalar potential which turns out to be minus some relative Fisher information, see Proposition \ref{res-12}. As the entropic interpolation, i.e.\,the marginal flow of the solution of the Schrödinger problem, solves the least action problem
\begin{align*}
\inf_{}\{ \mathsf{A}( \mu); \mu: \mu(0)= \mu_0, \mu(T)= \mu_T\},
\end{align*}
this suggests that it also 
solves some Newton equation in the Wasserstein space. See \cite{vR11, Co18} for some progresses in this direction.

It is worth mentioning that similar considerations apply to large deviation functionals of mean-field interacting particles (as opposed to non-interacting particle systems leading to the relative entropy $H(P|R)$), as for example in  \cite{BCGL19}. It brings us with a new interpretation in terms of Wasserstein geometry of the celebrated contributions of Dawson and Gärtner on the large deviations of mean-field particle systems \cite{DG87,DG89}.

Although this article focuses on time reversal, in order to clarify our motivation for studying  time reversal thirty-five years after it was well understood, we give some details about these considerations at Section \ref{sec-cod}, where Proposition \ref{res-12} is the rigorous statement of \eqref{eq-35a} and \eqref{eq-35b}.

\subsection*{Time reversal formula for a diffusion process}
General time reversal formulas for diffusion processes are well known since the 80's. 
Consider a diffusion process $Y$ in $\ZZ$ satisfying
\begin{align*}
dY_t= b_t(Y_t)\, dt+ \sigma_t(Y_t)\,dB_t,\quad 0\le t\le T,
\end{align*}
with $B$ a Brownian motion, $b$ a drift vector field and $ \sigma$ a matrix  field associated to the diffusion field $\aa:= \sigma \sigma ^{ \mathsf{t}},$ ($ \sigma ^{ \mathsf{t}}$ is the transposed of $ \sigma.$)  Assuming that the law  of $Y_t$ is absolutely continuous at each time $t$, under various hypotheses on $b$ and $\aa$, one can prove that the time-reversed process $Y^*$ is again a diffusion process with diffusion matrix field $\aa^*_t=\aa_{T-t}$ and drift field 
\begin{align}\label{eq-37}
b^*_t(y)= -b_{T-t}(y)+  \nabla\scal (\mu_{T-t} \aa_{T-t} )(y)/\mu_{T-t} (y),
\end{align}
where $\mu_t$ is the density of the law of $Y_t$ with respect to Lebesgue measure. This is not a straightforward result because a reversed semimartingale might not be a semimartingale anymore, see  \cite{W82}.

For this identity to hold,  it is assumed in \cite{HP86,MNS89} that   $b$ is locally Lipschitz (for a Sobolev-type  relaxation of this regularity property, see  \cite{RVW01}), and that  either $\aa$ is bounded away from zero or that the derivative $\nabla \aa$ in the sense of distribution is  controlled locally. Haussmann and Pardoux \cite{HP86} take a PDE approach, while Millet, Nualart and Sanz \cite{MNS89} rely on stochastic calculus of variations. The existence of an absolutely continuous density follows from a Hörmander type condition (PDE formulation in 
 \cite{HP86} and consequence of Malliavin calculus in \cite{MNS89}).
 
Föllmer's approach significantly departs from these strategies. Under the simplifying hypothesis that $\aa$ is the identity matrix, it is assumed in  \cite{Foe86} that the law $P$ of $Y$  has a finite entropy
\begin{align}\label{eq-36}
H(P|R)< \infty,
\end{align}
with respect to the law $R$ of a Brownian motion with some given initial probability distribution. In particular, the drift field $b$ of $P$  satisfies 
$
\IiZ |b_t(y)|^2\,\mu_t(y)dtdy < \infty
$
and  might be singular, rather than locally Lipschitz as required in  \cite{HP86,MNS89}. As a consequence of this finite entropy assumption, Föllmer proves  the time reversal formula
\begin{align}\label{eq-38b}
b^*_t(y)= -b_{T-t}(y)+  \nabla \log \mu _{ T-t}(y)
\end{align}
(recall $\aa=\Id$) where  the derivative is in the sense of distributions, without invoking any already known  result about the regularity of $\mu$.

With entropic optimal transport in mind, the hypothesis \eqref{eq-36} is mandatory. This rules out the Lipschitz regularity of $b$ which is required in ``non-Föllmerian'' approaches. Therefore, developing the entropic approach to time reversal is a necessary step of the research program attached to \emph{entropic optimal transport}. 
Following a previous unpublished work  by Cattiaux and Petit  \cite{CP01}, the present article fills this gap, keeping  the powerful guideline of Föllmer's proof based on stochastic derivatives and entropic arguments.

\subsection*{Main results of the article} Our main results are the IbP formula for the carré du champ of a general Markov process and the time reversal formula for a diffusion process.
\subsubsection*{IbP formula for the carré du champ of a Markov process}
Its expression is
\begin{align}\label{eq-24-intro}
E_P\Big((\LLf_t^Pu+\LLb_t^Pu)[X_t]v(X_t) + \Gaf_t^P(u,v)[X_t] \Big)
=0,
\end{align}
where $\LLf^P$, $\LLb^P$ are the forward and backward extended generators of the Markov measure $P,$ and $\Gaf^P$ is its forward extended carré du champ. See Section \ref{sec-sd} for more detail about these notions.  This IbP formula is valid for a sufficiently large class of  regular functions $u$ and $v.$
\\
No entropic argument is used to prove this result whose precise statement is given at Theorem \ref{res-02}. On the contrary, the main technical problem we face is to show that this IbP is valid under minimal regularity assumptions on $P$ to be able to apply it to general Markov measures typically satisfying a finite entropy condition. 
\\
The reason for calculating with extended generators is twofold:
\begin{enumerate}
\item
Unlike semigroup generators (which are associated to topological function spaces), extended generators are low-sophisticated objects which are tailor-made for \emph{martingale problems}: the relevant notion we work with in this article. This allows  us to consider lowly regular path measures $P$.

\item
As already noticed by Nelson in \cite{Nel67}, one can view Markov generators as stochastic derivatives, see Appendix \ref{sec-leo}. This natural idea  permits to perform computations along trajectories, using \emph{stochastic calculus} to obtain expressions for the generators and carré du champ operators. Our main technical result proved in this spirit is Lemma \ref{res-01}. It  is the keystone of the proof of the IbP formula.
\end{enumerate}

A discrete-time version of \eqref{eq-24-intro} was proposed by Feynman  \cite[Eq.\,(7-45)]{FH65} to derive Heisenberg's uncertainty principle with path integrals\footnote{We  thank Jean-Claude Zambrini for having brought this to our attention.}.

\subsubsection*{Time reversal formula for a diffusion process}

The law $P$ of the above process $Y$ solves the martingale problem 
\begin{align*}
P\in\MP(b,\aa)
\end{align*}
meaning that for any $u\in\CcZ$, the process $u(X_t)-\int_0^t \LLf_su(X_s)\,ds$ is a local $P$-martingale, where the forward generator $\LLf$   is defined  by
\begin{align*}
\LLf_tu(x)=b(t,x)\cdot \nabla u(x)+ \Delta_{\aa_t} u(x)/2,
\qquad (t,x)\in\iZ,
\end{align*}
with $\Delta_\aa:= \sum _{1\le i,j\le n} \aa_{ij} \partial^2_{ij}.$ One also writes \[P\in\MP(\mu,b,\aa)\] to specify the initial marginal measure $P_0=\mu$ if necessary. 
\\
The Markov generator of a \emph{Kolmogorov diffusion} with potential $U$ is
\begin{align*}
Au
 = \left( -\aa\nabla U\scal \nabla u+\nabla\scal(\aa\nabla u)\right)/2,
\end{align*}
where $\aa$ is a field on $\ZZ$ (not depending on time) with values in the set $S_+$ of all symmetric positive matrices and $U$ is a differentiable numerical function. The equilibrium measures of this dynamics are proportional to 
\begin{align*}
\mm(dx)=e ^{ -U(x)}\,dx.
\end{align*}
Expanding the divergence term, we see that the drift field of the generator  is 
\begin{align*}
\vam := (\nabla\scal \aa-\aa\nabla U)/2.
\end{align*}
\begin{hypotheses}\label{ass-01}\ 
\begin{enumerate}[(i)]
\item
$U\in C^1(\ZZ)$, $\aa$ is invertible and in $C ^1(\ZZ,S_+)$,
\item
for  some $K\ge 0,$\quad  $x\scal \vv ^{ \aa,\mm}(x)+\tr \aa(x)\le K(1+|x|^2)$ for all $x\in\ZZ$.
\end{enumerate}
\end{hypotheses}
 It is a standard result that under these hypotheses, 
the martingale problem $\MP(\mm,\vv ^{ \aa,\mm},\aa )$ admits a unique solution denoted by
\begin{align}\label{eq-09}
R \in\MP(\mm,\vv ^{ \aa,\mm},\aa ),
\end{align}  
which is  $\mm$-reversible. This implies in particular that $R^*=R.$

\begin{theorem}[Time reversal formula]\label{res-08-intro}
Under the Hypotheses \ref{ass-01} on $R$ given at \eqref{eq-09}, let $P\in\PO$ be Markov and such that \[H(P|R)< \infty.\]
Then, for all $t$ the density $ \mu_t:=dP_t/d\Leb$ exists and the time reversal $P^*$ of $P$ is a solution of the martingale problem 
\begin{align*}
P^*\in \MP(b^*,\aa)
\end{align*} 
with 
\begin{align}\label{eq-13d-intro}
b^*_t(x)	
	=-b_{T-t}(x)+   \nabla\scal( \mu_{T-t} \aa)(x)/ \mu_{T-t}(x),
	\quad dtP_t(dx)\ae
\end{align}
where the divergence is in the sense of distributions.  
\\
This is an extension of  \eqref{eq-37}  to a low regularity setting which is made precise as follows.
\\
For almost every $t$ the density $\rho_t:=dP_t/d\mm$ admits a distributional spatial derivative $\nabla \rho_t$  satisfying
\begin{align}\label{eq-15b-intro}
\IiZ |\nabla \log \rho_t|^2_\aa \, dP_tdt< \infty.
\end{align}
and \eqref{eq-13d-intro} is equivalent to
\begin{align} \label{eq-fundTR}
(b_t+b^*_{T-t})/2-\vam=\aa \nabla\log \sqrt{\rho_t},
 \quad dtdP_t\ae
\end{align}
Furthermore, $P^*$ is the unique solution of $\MP(P_T,b^*,\aa)$ among the set of  all $Q\in\PO$ such that $H(Q|R) < \infty.$
\end{theorem}

This theorem is a restatement of Theorem \ref{res-08} which is stated  in terms of stochastic velocities, especially the fundamental identity \eqref{eq-fundTR} which is synthetically expressed in terms of the \emph{osmotic momentum} $\bo{P|R}$ at \eqref{eq-14c}.

\subsection*{An extension of Theorem \ref{res-08-intro}}

Note that unlike \cite{HP86,MNS89},  it is assumed in Theorem  \ref{res-08-intro}  that the diffusion matrix field $\aa$ does not depend on $t$. However, our method allows to extend the results of \cite{HP86,MNS89} to a finite entropy setting. Indeed, the method of proof of the present article is perturbative: \emph{if one knows the time-reversal formula for some possibly unbounded reference path measure $R$, then a time-reversal also holds for any probability measure $P$ such that $H(P|R) < \infty.$}

A careful inspection of the proof of Theorem \ref{res-08} (a.k.a.\,Theorem \ref{res-08-intro}) shows that it extends to the case where the reference measure $R$ might not be reversible.

\begin{theorem}[Time-reversal formula, again] \label{res-08-extension} Let us assume that the possibly unbounded reference measure $R$ and  its  time reversal $R^*$ both solve uniquely their respective  martingale problems $\MP(R_0,b^R,\aa)$ and $\MP(R^*_0=R_T,b ^{ R^*},\aa^*)$ in the sense of \eqref{eq-H17a},  where $b^R, b ^{ R^*}$ are locally bounded fields  and $\aa$ is continuous on $\iZ$. 
The following assertions are verified.
\begin{enumerate}[(a)]
\item
For  all $0\le t\le T,$ we have   $\aa^*_t=\aa _{ T-t}.$  
\item
Assume also  that for all $0< t< T$    the time marginal $R_t$ is absolutely continuous with respect to Lebesgue measure, and that the path probability  measure satisfies $H(P|R)< \infty$ again. \\ Then, $P$ and $P^*$  uniquely solve  $\MP(b^P,\aa) $ and $\MP(b ^{ P^*},\aa^*)$ respectively, in the sense of  \eqref{eq-H17a}. The identity   \eqref{eq-fundTR} becomes
\begin{align}
(b_t^P+b ^{ P^*}_{T-t})/2- (b_t^R+b ^{ R^*}_{T-t})/2=\aa_t \nabla\log \sqrt{\rho_t}, \quad dtdP_t\ae,
\end{align}
where $ \rho_t= d P_t/dR_t$ and the gradient is in the  sense of distribution, and \eqref{eq-15b-intro}  still holds:
\begin{align*}
\IiZ |\nabla \log \rho_t|^2 _{ \aa_t} \, dP_tdt< \infty.
\end{align*}
\end{enumerate}
 \end{theorem}

\begin{proof}
Statement (a) is a consequence of \eqref{eq-45} at Lemma \ref{res-25}-(b) applied to $R$, whose assumptions are satisfied by Theorem \ref{res-02}-(b).
\\
The proof of item (b)  is similar to the proof of Theorem \ref{res-08}, almost verbatim. The uniqueness of the solution to the martingale problems for $R$ and $R^*$ is necessary for invoking Girsanov's theory at Proposition \ref{res-05}.  Finally, the local boundedness of the semimartingale characteristics of $R$ and $R^*$ implies  the boundedness of  $\LLf^R u$ and $\LLb^R u$ for any $u\in\CcZ.$  This enters the proof of Lemma \ref{res-06} in an essential manner.
\end{proof}

\begin{remark} \label{rem-UG}\
A typical hypothesis for a path measure $Q$ to be the unique solution of its martingale problem $\MP(\aa,b^Q)$ in the sense of \eqref{eq-H17a} is that $\aa= \sigma\,\sigma^*$ with $\sigma$ and $b^Q$ locally Lipschitz in space and time.
\end{remark}

In particular, with  $R$ satisfying the regularity hypotheses of the main results of \cite{HP86,MNS89}, we see that Theorem \ref{res-08-extension} extends the time reversal formula \eqref{eq-37} to the  wider class of all path measures $P$ such that $H(P|R) < \infty.$

\subsection*{Literature about time reversal of Markov processes}

The first investigations in the theory of time reversal of Markov processes date 
back to 1936 with a pair of articles \cite{Kol36,Kol37} by Kolmogorov providing sufficient conditions for a Markov chain or a diffusion process to be reversible. Then, in 1958 time reversal of Markov processes was used by Hunt  \cite{H58} in his study of potential theory. During the same year, Nelson published an article \cite{N58} entitled ``The adjoint Markoff process''. Several papers  went on in the direction initiated by Hunt: \cite{N64, KW66,D85} (to cite a few of them). All these articles deal with  \emph{stationary} Markov processes and their results are expressed in terms of transition probabilities  rather than semimartingale characterics, which is quite natural in the framework of potential theory.  The above mentioned articles \cite{Foe86,HP86,MNS89} and \cite{Par86} are the first ones where the expression  of  semimartingale characteristics of a time reversed process are obtained rigorously, see \eqref{eq-37}. They are restricted to a diffusion setting.
\\
 We also mention the  article \cite{KST18} by Karatzas, Schachermayer and Tschiderer both for its  well written appendix section  on time reversal of diffusion processes and its results connecting  deterministic optimal transport and diffusion processes, where time reversal plays a crucial role. The recent article \cite{KMS20} by Karatzas, Maas and Schachermayer also makes use of time reversal in the context of Markov chains.
Recently,  Izydorczyk, Oudjane, Russo and  Tessitore \cite{IOR21,IORT21} used time reversal of diffusion processes  to prove the well-posedness of some backward Fokker-Planck equations, and to design  efficient algorithms  solving some Hamilton-Jacobi-Bellman equations with terminal conditions.

\subsubsection*{Nelson's contribution}

While investigating large deviations of the empirical measure of weakly interacting Brownian particles  as in \cite{Sch31,Sch32} or  \cite{DG87}, Föllmer   established the time reversal formula \eqref{eq-38b} using entropic arguments, among which the identity \eqref{eq-38} is decisive.  At the same period,  Zambrini obtained in \cite{Zam86} a time-symmetric description of the backward and forward drifts of the solution to the  Schrödinger problem \eqref{eq-Sch-pb}. These two authors used in a crucial manner the notion of stochastic derivatives introduced by Nelson in  1967 in \cite{Nel67}. Time reversal is at the core  of Nelson's theory of Brownian motion. Indeed, his expression of the osmotic velocity (a notion introduced by him, guided by the seminal article \cite{Ein05} by Einstein on the Brownian motion)  in terms of the density of the process is nothing but the time reversal formula. He proves it in an informal manner, i.e.\  assuming that all the derivatives exist in a classical sense, using PDEs, namely Fokker-Planck equations in both directions of time, also called forward and backward Kolmogorov's equations after   \cite{Kol37}.  In the present article, stochastic derivatives also play a major role.

\subsubsection*{Back to the roots}

As a concluding remark about  the history of time reversal of Markov processes, it appears that the very starting point of this adventure is, again, the paper \cite{Sch31} by Schrödinger. Indeed, in the first paragraph of  \cite{Kol36},  Kolmogorov refers to \cite{Sch31} as his main motivation\footnote{We  thank Jean-Claude Zambrini for having brought this to our attention.}.

\subsection*{Notation}

The set of all probability measures on a measurable set $A$ is denoted by $ \mathrm{P}(A)$ and the set of all nonnegative $ \sigma$-finite measures on $A$ is  $ \mathrm{M}(A).$
The push-forward of a measure $ \qq\in  \mathrm{M}(A)$ by the measurable map $f:A\to B$ is   $f\pf\qq(\sbt):=\qq(f\in \sbt)\in \mathrm{M}(B).$ 

\subsubsection*{Relative entropy}

The relative entropy of $\pp\in \mathrm{P}(A)$ with respect to the reference  measure $\rr\in \mathrm{M}(A)$ is 
\begin{align*}
H(\pp|\rr):=\int_A \log(d\pp/d\rr)\,d\pp \in( - \infty, \infty]
\end{align*}
if $\pp$ is absolutely continuous with respect to $\rr$ ($\pp\ll\rr$) and $\int_A \log_-(d\pp/d\rr)\,d\pp< \infty$, { and $H(\pp|\rr)=+ \infty$ otherwise}. If $\rr\in \mathrm{P}(A)$ is a probability measure, then $H(\pp|\rr) \in[0, \infty].$ See Section \ref{sec-um} for details.

\subsubsection*{Path measures}
The configuration space is a Polish space $\XX$ equipped with its Borel $ \sigma$-field.
The path space is the set { $\OO:=D(\ii,\XX)$} of all $\XX$-valued \cadlag\  trajectories on the time index set $\ii,$ and the canonical process  $(X_t) _{ 0\le t\le T}$ is defined by $X_t( \omega)= \omega_t$ for any $0\le t\le T$ and any path $ \omega=( \omega_s) _{ 0\le s\le T}\in\OO.$ It is equipped  with the canonical $ \sigma$-field $ \sigma(X _{ \ii})$ and  the canonical filtration $ \Big(\sigma( X _{ [0,t]}); 0\le t\le T\Big)$ where for  any subset $ \mathcal{T}\subset \ii$, $X _{ \mathcal{T}}:=(X_t, t\in \mathcal{T})$ and  $ \sigma(X _{ \mathcal{T}})$ is the $ \sigma$-field generated by the collection of maps $(X_t, t\in \mathcal{T})$.  
\\
The càdlàg setting is necessary at Section \ref{sec-df} for the abstract IbP formula and Section \ref{sec-rw} where random walks are investigated. At Sections \ref{sec-diff1},  \ref{sec-diff2} and \ref{sec-cod}, diffusion processes are time-reversed and the path space is the set  $\OO=C(\ii,\ZZ)$ of all continuous trajectories.
\\
 We call any positive measure  $Q\in\MO$ on $\OO$ a path measure. For any $\mathcal{T}\subset\ii,$ we denote  $Q_\mathcal{T}=(X_\mathcal{T})\pf Q.$ In particular, for any
 $0\le r\le s\le T,$ $X_{[r,s]}=(X_t)_{r\le t\le s}$, $Q_{[r,s]}=(X_{[r,s]})\pf Q$, and $Q_t=(X_t)\pf Q\in\MX$ denotes the law of the position $X_t$ at time $t$. If $Q\in\PO$ is a probability measure, then $Q_t\in\PX$.
\\
The time-space canonical process is
$$
\Xb_t:=(t,X_t)\in\iX,
$$
and for  any function $u:\iX\to\RR$,  we denote  $u(\Xb): (t, \omega)\mapsto u(t, \omega_t).
$
We also denote
\begin{align*}
\Qb  (dtd \omega)&:=dtQ (d \omega),\qquad dtd \omega\subset\ii\times \OO,\\
\qb (dtd x)&:=dtQ_t (d x),\qquad dtd x\subset\iX.
\end{align*}

\section{Stochastic derivatives} \label{sec-sd}

Let us recall the definitions of Markov measures, extended generators and stochastic derivatives. The precise definitions of these notions together with some useful related technical results are recalled at the appendix Section \ref{sec-leo}. Stochastic derivatives were introduced by Nelson in 1967 \cite{Nel67}.

\subsection*{Conditionable path measure}

  A path measure $Q$ such that  $Q_t$ is  $ \sigma$-finite for all $t$ is called a conditionable path measure. This notion is necessary to define properly the conditional expectations $E_Q(\sbt\mid X _{t }),$ $E_Q(\sbt\mid X _{[0,t] })$ and $E_Q(\sbt\mid X _{[t,T] }),$ for any $t$, see \cite{Leo12b}. If $Q$ has a finite mass, then it is automatically conditionable.

\subsection*{Extended forward generator}

 Let $Q$ be a conditionable measure. A measurable function $u$ on $\iX$ is said to be in
the domain of the extended forward generator of $Q$ if there exists a real-valued 
process $\LLf^Q u(t,X _{ [0,t]})$  which is adapted with respect to the forward filtration such that
 $\Iii|\LLf^Q u(t,X _{ [0,t]})|\,dt<\infty,$ $Q\ae$ and the process $$M^u_t:=u(\Xb_t)-u(\Xb_0)-\II0t \LLf^Q u(s,X _{ [0,s]})\,ds,\quad 0\le t\le T,$$ is a local $Q$-martingale.
We say that $\LLf^Q $ is  the extended forward generator of $Q.$ Its domain  is denoted by $\dom
\LLf^Q .$ 
Otherwise stated, we say that $Q$ solves the \emph{martingale problem} 
with generator $\LLf$ and domain $\UU,$ 
if $\UU\subset \dom\LLf^Q$ and for any $u\in \UU,$ $\LLf^Qu=\LLf u.$

\subsection*{Stochastic forward derivative}

Nelson's definition \cite{Nel67} of the stochastic forward derivative is the following. 
For any   conditionable measure $Q$ and any  measurable real function $u$  on $\iX$ such that $E_Q|u(\Xb_s)|<\infty$ for all $0\le s\le T,$ 
we say
that $u$ admits a stochastic forward derivative  under $Q$ at time
$t\in[0,T)$ if   the following limit
\begin{align}\label{eq-100}
 \Lf^Qu(t,X _{ [0,t]}):=\Limh E_Q\left(\frac1h [u(\Xb_{t+h})-u(\Xb_t)]
    \mid  X _{ [0,t]}  \right)
\end{align}
exists in $L^1(\Qb).$
In this case,  $\Lf^Qu(t,\sbt)$ is called  the stochastic forward derivative of $u$ at time $t$.

\subsection*{Extended generators and stochastic derivatives are essentially the same}

It is the content of Proposition \ref{res-10}.
  If  $u$ is in $\dom\LLf^Q$ and satisfies
$E_Q\Iii
 \big|\LLf^Q u(t,X _{ [0,t]})\big|\,dt<\infty$, one can compute  $\LLf^Q u$  using the stochastic derivative:
\begin{align*}
\LLf^Q u=\Lf^Qu,\quad \Qb\ae
\end{align*}
Beware of the notation: calligraphic $\LL$ refers to the martingale problem $\MP(\LL),$ while the roman font $L$ refers to the stochastic derivative \eqref{eq-100} which provides us with a mean of calculating $\LL$ via \eqref{eq-100} using  stochastic calculus.

\subsection*{Reversing time}

Let $Q\in\MO$ be any path measure. Its time reversal is 
\begin{align*}
Q^*:=(X^*)\pf Q\in\MO,
\end{align*}
where 
\begin{align*}
\left\{ \begin{array}{ll}
X^*_t:= \Limh X _{ T-t+h},\quad & 0\le t< T,\\
 X^*_T:=X_0,& t=T,
\end{array}\right.
\end{align*}
is the reversed canonical process. We assume that $Q$ is such that 
\[
Q(X _{ T^-}\neq X_T)=0,
\] 
i.e.\ its sample paths are left-continuous at $t=T.$ This implies that the time reversal mapping $X^*$ is (almost surely) one-one on $\OO.$

As a notation, the $ \sigma$-field generated by $ X _{ [t^-,T]}$ is $ \sigma(X _{ [t^-,T]}):= \cap _{ h>0} \sigma(X _{ [t-h,T]})= \sigma(X _{ t^-})\vee \sigma(X _{ [t,T]}),$ and the predictable backward filtration is defined by: $( \sigma(X _{ [t^-,T]}); 0\le t\le T).$
 
 We introduce the backward extended generator and the backward stochastic derivative
\begin{align}\label{eq-F05}
\begin{split}
\LLb^Qu(t, X _{ [t^-,T]})&:=\LLf ^{ Q^*} u^*(t^*,X^* _{ [0,t^*]}),\\
\Lb^Qu(t, X _{ [t^-,T]})&:=\Lf ^{ Q^*} u^*(t^*,X^* _{ [0,t^*]}),
\end{split}
\end{align}
where $ u^*(t^*, x):=u(t, x),$  with $t^*:=T-t$,   and $ \LLf ^{ Q^*}$ and  $\Lf ^{ Q^*}$ stand respectively for the standard (forward) generator and derivative of $Q^*$. These definitions  match with  Definitions \ref{def-02b} and \ref{def-03b}. In particular, for any $t\in (0,T],$
\begin{align*}
 \Lb^Qu(t,X _{ [t^-,T]}):=\Limh E_Q\left(\frac1h [u(\Xb_{t-h})-u(\Xb_t)]
    \mid  X _{ [t^-,T]}  \right)
\end{align*}
if this limit exists in $L^1(Q).$  
Remark that the definition of $\Lb^Q$ is consistent with \eqref{eq-F05}.

The linear operators $\LLf, \Lf, \LLb$ and $\Lb$  are defined for any measurable function $u:\iX\to\RR$ such that the above expressions are meaningful where this \emph{meaningful} addresses the   problem of their domains, see the appendix section \ref{sec-leo}.
\\
As for the forward generator,  if  $u$ is in $\dom\LLb^Q$ and satisfies
$E_Q\Iii
 \big|\LLb^Q u(t, X _{ [t,T]})\big|\,dt<\infty$, then
 \begin{align*}
\LLb^Q u=\Lb^Qu,\quad \Qb\ae
\end{align*}

  \subsection*{Markov measure}

  A path measure $Q\in\MO$ is said to be Markov if it is conditionable and for any $0\le t\le T,$ $Q(X _{ [t,T]}\in\sbt\mid X _{ [0,t]})=Q(X _{ [t,T]}\in\sbt\mid X _t).$  It is known that $ Q^*$
 is also Markov and the stochastic derivatives and extended generators at time $t$ only depend of the present position $X_t$. Therefore it is possible to consider the sum and difference of the forward and backward generators: they remain functions of the present position.

\subsubsection*{Current and osmotic generators}
In restriction to  $\dom\LLf^Q\cap \dom\LLb^Q,$ we define the \emph{current} extended generator of $Q$ by
\begin{align*}
\LLc Q:= (\LLf^Q-\LLb^Q)/2.
\end{align*}
Similarly, the \emph{osmotic} extended generator of $Q$ is
\begin{align*}
\LLo Q:= (\LLf^Q+\LLb^Q)/2.
\end{align*}
The osmotic generator plays an important role in this article. This is the reason why our results about time reversal are restricted to Markov measures.

\section{Integration by parts formula} \label{sec-df}

The main technical result of this paper is the integration by parts formula stated  at Theorem \ref{res-02}.  This section is dedicated to its statement and its proof.

\subsection*{Carré du champ}

Let $Q$ be a path measure on $\OO.$
Its forward  carré du champ  is the forward-adapted process defined by
\begin{align*}
\Gaf^Q_t(u,v):= \LLf^Q_t(uv)-u\LLf^Q _tv-v\LLf^Q_t u,
\quad (u,v)\in\dom \Gaf^Q_t,\ 0\le t\le T,
\end{align*}
where 
$
\dom\Gaf^Q_t:= \left\{(u,v); \  u, v, uv \in\dom \LLf^Q_t \right\} .
$
\\
We introduce a class $\UU$ of functions on $\XX$ such that
\begin{align}\label{eq-101a}
 \mathcal{U}\subset \dom \LLf^Q_t\cap C_b(\XX)
\end{align}
for all $0\le t\le T$ and any path measure $Q$ of interest, where $C_b(\XX)$ is the space of all bounded continuous functions on $\XX.$  \emph{We assume that $\UU$ is an algebra,} i.e. 
\begin{align}\label{eq-101b}
u,v\in\UU\implies uv\in\UU.
\end{align}
In particular, 
\begin{align}\label{eq-43}
u,v\in\UU\implies (u,v)\in\dom \Gaf^Q _t.
\end{align}
We shall mainly consider functions in $\UU$ and make an intensive use of their carré du champ. In each setting, this algebra will be chosen rich enough to determine a Markov dynamics, i.e.\ to solve in a unique way some relevant  martingale problem. For instance, in the diffusion setting, $\UU=\CcZ $ is a good choice.

\begin{remark}\label{rem-04}
The requirement that {$ \mathcal{U}$ is an algebra} (it is necessary that $uv$ belongs to $\dom \LLf^Q$ to consider $\LLf^Q (uv)$), is  strong. 
\\
Indeed, let us say that  a semimartingale is \emph{nice} if its bounded variation part is \emph{absolutely continuous.} The product of two semimartingales is a semimartingale, but  the product of two {nice} semimartingales might not be  nice anymore. 
\\
However, this is true for instance when the semimartingales are adapted to a Brownian filtration because in this case any local martingale is represented as a stochastic integral with respect to a Brownian motion. In general, a martingale representation theorem is needed to verify the stability of the product of nice semimartingales.
\end{remark}

Similarly  the backward carré du champ  is the backward-adapted process defined by
\begin{align*}
\Gab^Q_t(u,v):= \LLb^Q_t(uv)-u\LLb^Q_t v-v\LLb^Q_t u,
\end{align*}
for any $0\le t\le T$ and $(u,v)\in\dom \Gab^Q_t.$ 
To emphasize the fact that $\Gaf^Q(u,v)$ and $\Gab^Q(u,v)$ are processes rather than functions, we often write 
\begin{align*}
\Gaf^Q_t(u,v)=\Gaf^Q_t(u,v)(X _{ [0,t]})=\Gaf^Q_t(u,v)(X),\\
\Gab^Q_t(u,v)=\Gab^Q_t(u,v)(X _{ [t,T]})=\Gab^Q_t(u,v)(X).
\end{align*}
The quadratic covariation $[u(X),v(X)]$ is a $Q$-semimartingale.
We denote by $\langle u(X),v(X)\rangle^Q$ its bounded variation part,  i.e.
\begin{align}\label{eq-41}
d[u(X),v(X)]_t=d\langle u(X),v(X)\rangle^Q_t+ dM_t ^{Q, [u,v]},
\qquad \Qb\ae
\end{align}
where, here and below, $M^{Q,\sbt}$  stands for any forward local  $Q$-martingale. As next lemma indicates, we are interested in situations where the bounded variation process $\langle u(X),v(X)\rangle^Q$ is predictable (as a continuous process). Therefore, in the remainder of the article $\langle u(X),v(X)\rangle^Q$ is the usual sharp bracket (sometimes called conditional quadratic variation) of stochastic process theory.

\begin{lemma}\label{res-13}
Let $\UU$ satisfy the hypotheses \eqref{eq-101a} and \eqref{eq-101b}.
\begin{enumerate}[(a)]
\item
For any $u,v\in\UU,$ the process $ \langle u(X),v(X)\rangle^Q $ is absolutely continuous $Q\ae$ and 
\begin{align*}
 d\langle u(X),v(X)\rangle^Q_t=\Gaf^Q_t(u,v)(X _{ [0,t]})\,dt, \qquad \Qb\ae
\end{align*}

\item
For any $u,v\in\UU,$ the process $ \langle u(X),v(X)\rangle ^{ Q^*} $ is absolutely continuous $Q^*\ae$ and 
\begin{align*}
 d \langle u(X),v(X)\rangle ^{Q^*} _{ |T-t}(X^*)=\Gab^Q_t(u,v)(X _{ [t,T]})\,dt, \qquad \Qb\ae
\end{align*}
\end{enumerate}
\end{lemma}

\begin{proof}
\boulette{(a)}
As a definition of the forward generator
\begin{align*}
&du(X)_t=\LLf^Q_tu(X)\,dt + dM^u_t,
\qquad dv(X)_t=\LLf^Q_t v(X)\, dt+dM^v_t,\\
&d(uv)(X)_t=\LLf^Q_t(uv)(X)\,dt+dM ^{ uv}_t,
\end{align*}
and  applying Itô's formula in the forward sense of time \begin{align*}
d(uv)(X)_t&=u(X_t)dv(X)_t+v(X_t)du(X)_t+d[u(X),v(X)]_t\\
	&=u(X_t)dv(X)_t+v(X_t)du(X)_t+d\langle u(X),v(X)\rangle_t+dM_t ^{Q, [u,v]}\\
	&=[u(X_t)\LLf^Q_tv(X)+v(X_t)\LLf^Q _tu(X)]\,dt+d\langle u(X),v(X)\rangle_t\\
	&\hskip 4cm +u(X_t)dM^v_t+v(X_t)dM^u_t+dM_t ^{Q, [u,v]}.
\end{align*}
The Doob-Meyer decomposition theorem allows us to identify the bounded variation and martingale parts of $uv(X)$, leading us to
\begin{align*}
\LLf^Q_t(uv)(X)\, dt
	= [u(X_t)\LLf^Q_tv(X)+v(X_t)\LLf^Q _tu(X)]\,dt+d\langle u(X),v(X)\rangle^Q_t,
	\qquad \Qb\ae
\end{align*}
which gives the announced result.
\Boulette{(b)}
Analogous, with $Q^*$ instead of $Q.$
\end{proof}

Remark that the main hypothesis of this lemma is \eqref{eq-101b}: $u,v\in\UU$, and its consequence \eqref{eq-43}.

Let us  prepare some notation for  next  Lemma \ref{res-01} which is the main technical result of this section.
We 
introduce the class of functions
\begin{equation}\label{eq-42b}
\UQ_2 := \Big\{ u\in \UU; \LLf^Qu(X)\in L^2(\Qb),\ \Gaf^Q (u)(X)\in L^1(\Qb)  \Big\}.
\end{equation}
If $Q$ is Markov, $\Gaf^Q_t(u,v)(X)=\Gaf^Q_t(u,v)[X_t] $ only depends on the current position $X_t$, and we denote  $$(t,x)\mapsto \Gaf^Q_t(u,v)[x]:=E_Q(\Gaf^Q_t(u,v)(X)\mid X_t=x).$$  
  Consider the following convolution kernels
\begin{align*}
k^h:= h ^{ -1} \1 _{ [-h,0]},
\qquad
k ^{ -h}:= h ^{ -1}\1 _{ [0,h]},
\end{align*}
with $h>0$. 
Let $f:\ii\to \RR$ be any absolutely continuous function  with derivative $\dot f.$ The following expressions will be used during the proof of next lemma:
\begin{align}\label{eq-22}
\begin{split}
&h ^{ -1}[f(t+h)-f(t)]= h ^{ -1}\int _{ [t,t+h]}\dot f(r)\,dr= k^h\ast \dot f(t),\quad 0\le t\le T-h,\\
&h ^{ -1}[f(t)-f(t-h)]= h ^{ -1}\int _{ [t-h,t]}\dot f(r)\,dr= k ^{ -h}\ast \dot f(t),\quad h\le t\le T.
\end{split}
\end{align}

\begin{lemma}\label{res-01}
Let $Q$ be any  path measure and take any $u,v$ in the class  $\UQ_2.$
\begin{enumerate}
\item[(a)]
The following limit holds
 \begin{multline*}
\Limh E_Q\int_0 ^{ T-h} \Big|E_Q \big[h ^{ -1}\{u( X_{t+h})-u( X_t)\}\{v( X_{t+h})-v( X_t)\}\mid X _{ [0,t]}\big]\\
	-   \Gaf^Q_t(u,v)(X)\Big|\,dt =0.
\end{multline*}

\item
[(b)]
If in addition $Q$ is Markov and    $(t,x)\mapsto \Gaf^Q_t(u,v)[x] $ is continuous, then 
\begin{equation} \label{eq-44}
\begin{split}
\Limh E_Q\int_h ^{ T} \Big|E_Q \big[h ^{ -1}\{u( X_{t})-u( X _{ t-h})\}\{v( X _{ t})-v&( X _{ t-h})\}\mid X _{ t-h}\big]\\
	&-   \Gaf^Q _{ t}(u,v)[X_t]\Big|\,dt =0.
\end{split}
\end{equation} 
\end{enumerate}
\end{lemma}

\begin{proof}
\boulette{(a)}
Let us start with a remark about our assumptions. The (a priori local) martingale $M^u_t=u(X_t)-u(X_0)-\int_0^t\LLf^Q_su( X)\,ds$, is a square integrable martingale because
\begin{align}\label{eq-21}
E_Q \sup _{ 0\le t\le T} |M^u_t|^2
	\le C_2 E_Q {[ u(X)]^Q_T}
	= C_2 {E_Q \langle u(X)\rangle ^Q_T}
	= C_2 T {E_Q  \Gaf^Q_T(u)(X)}< \infty,
\end{align}
where the first inequality is Doob's maximal inequality with $C_2=4$, and the rest follows from   the assumptions $\Gaf^Q (u)(X)\in L^1(\Qb)$ and  Lemma \ref{res-13}.
For each $0\le t\le T-h$ with $0<h\le T,$
\begin{align*}
[u( X_{ t+h})-u( X_t)]&[v( X _{ t+h})-v( X_t)]\\
	= &\Big[ \int_t ^{ t+h} dM^u_s
		+ \int_t ^{ t+h} \LLf^Q_su(X)\,ds\Big]
		 \Big[ \int_t ^{ t+h} dM^v_s
		+ \int_t ^{ t+h} \LLf^Q_sv(X)\,ds\Big]\\
	= & A_t^h+B_t^h+C_t^h+D_t^h,
	\qquad Q\ae,
\end{align*}
where
\begin{align*}
A_t^h&= \int_t ^{ t+h} dM^u_s\ \int_t ^{ t+h} dM^v_s,\hskip 2,5cm
B_t^h=  \int_t ^{ t+h} \LLf^Q_su(X)\,ds\  \int_t ^{ t+h} dM^v_s,\\
C_t^h&=  \int_t ^{ t+h} \LLf^Q _sv(X)\,ds\  \int_t ^{ t+h} dM^u_s,\hskip 1,1cm 
D_t^h=  \int_t ^{ t+h}\LLf^Q_su(X)\,ds\  \int_t ^{ t+h} \LLf^Q_sv(X)\,ds.
\end{align*}
Let us control $A^h_t.$ Denoting $N^u_{t,s}:=M^u_s-M^u_t$ and 
$N^v_{t,s}:=M^v_s-M^v_t,$
\begin{align*}
A_t^h&=
	\int_t ^{ t+h} d(N^u_{t,s}N^v_{t,s})\\
&= \int_t ^{ t+h} N^u_{t,s}  dM^v_s
		+ \int_t ^{ t+h} N^v_{t,s}  dM^u_s
		+ \int_t ^{ t+h} dM ^{Q, [u,v]}_s
		+ \int_t ^{ t+h} d \langle M^u,M^v\rangle ^Q _s,
\end{align*}
where $M ^{ Q,[u,v]}$ is the martingale part of the semimartingale $[u(X),v(X)]$, see \eqref{eq-41}.
With Lemma \ref{res-13}, we obtain
\begin{align}\label{eq-02}
h ^{ -1}E_Q (A^h_t\mid X _{ [0,t]})= h ^{ -1}\int_t ^{ t+h} E_Q[\Gaf^Q_s(u,v)(X _{ [0,s]})\mid X _{ [0,t]}]\,ds.
\end{align}
Remark that under our integrability assumptions, the stochastic integrals  $\int_t ^{ t+h} dM^u_s$, $\int_t ^{ t+h} dM^v_s,$ $\int_t ^{ t+h} N^u_{t,s}  dM^v_s$ and $\int_t ^{ t+h} N^v_{t,s}  dM^u_s$  are integrable $Q$-martingales. The first ones because of  \eqref{eq-21}, and the last ones by Burkholder-Davis-Gundy inequality:
\begin{multline*}
E_Q \sup _{ 0\le t\le T} \Big|\int_0^t M^u_s\,dM^v_s\Big|
	\le C_1 E_Q\left[ \left( \int_0^T |M^u_t|^2d[M^v]_t \right) ^{ 1/2}\right]
	\le C_1 E_Q \left(\sup _{ 0\le t\le T}  |M^u_t| [M^v]_T^{ 1/2}\right) \\
	\le C_1  \sqrt{E_Q \sup _{ 0\le t\le T}  |M^u_t| ^2} \sqrt{\ E_Q [M^v]_T}
	\le  C_1 C_2 ^{ 1/2}   \sqrt{E_Q \langle u(X)\rangle_T} \sqrt{\ E_Q \langle v(X)\rangle_T}
	< \infty,
\end{multline*}
with $C_1$ a universal constant and where we used \eqref{eq-21} when $C_2$ appears.
 By Burkholder-Davis-Gundy inequality again, we also have $M ^{ Q,[u,v]}\in L^1(Q).$ These considerations  justify the cancellation of the expectations of the martingale terms.

The remaining terms $B^h, C^h$ and $D^h$ are controlled using our integrability assumptions and Cauchy-Schwarz inequality.
Let us start with $B^h$:
\begin{align*}
&\ \Big(E_Q \int_0 ^{ T-h}|B_t^h|\,dt\Big)^2\\
	&\le  E_Q \int_0 ^{ T-h}\big(\int_t ^{ t+h} \LLf^Q_su(X)\,ds\big)^2\,dt\ \ 
	E_Q \int_0 ^{ T-h}\big(\int_t ^{ t+h} dM^v_s\big)^2\,dt
	\\
	&\le  E_Q \int_0 ^{ T-h}\big(\int_t ^{ t+h} \LLf^Q_su(X)\,ds\big)^2\,dt\ \ 
	E_Q\int_0 ^{ T-h}\int_t ^{ t+h} \Gaf^Q_s(v)(X)\,dsdt\\
	&\le o(h^2)\ 
	E_Q \int_0 ^{ T-h}  k^h\ast (\LLf^Q u)^2(t,X _{ [0,t]})\,dt\\
	&= o(h^2) \ 
	\Big(E_Q \int_0 ^{ T} (\LLf^Q_tu)^2(X)\,dt+o_{h\to 0^+}(1)\Big),
\end{align*}
where the third inequality follows from Lebesgue's dominated convergence theorem under the assumption that $\Gaf^Q(v)(X)\in L^1(\Qb)$, and  use  we took $k^h:= h ^{ -1} \1 _{ [-h,0]}$ as our convolution kernel, see \eqref{eq-22}. The last identity is a consequence of  Lemma \ref{res-09} under the assumption  $\LLf^Qu(X)\in L^2(\Qb)$. This gives 
\begin{align*}
E_Q \int_0 ^{ T-h}h ^{ -1}|B_t^h|\,dt
	\le  o_{h\to 0^+}(1) \|\LLf^Qu(X)\| _{ L^2(\Qb )}+o_{h\to 0^+}(1)
\end{align*}
and similarly
\begin{align*}
E_Q \int_0 ^{ T-h}h ^{ -1}|C_t^h|\,dt
	\le o_{h\to 0^+}(1) \|\LLf^Qv(X)\| _{ L^2(\Qb )}+o_{h\to 0^+}(1).
\end{align*}
The control of $D^h$  is analogous:
\begin{align*}
&\ \Big(E_Q \int_0 ^{ T-h}|D_t^h|\,dt\Big)^2\\
	&\le  E_Q \int_0 ^{ T-h}\big(\int_t ^{ t+h} \LLf^Q_su(X)\,ds\big)^2\,dt\ 
	E_Q \int_0 ^{ T-h}\big(\int_t ^{ t+h} \LLf^Q_sv(X)\,ds\big)^2\,dt
	\\
	&\le h^4 E_Q \int_0 ^{ T-h}k^h\ast (\LLf^Qu)^2(t,X _{ [0,t]})\,dt\ 
	E_Q \int_0 ^{ T-h}k^h\ast (\LLf^Q_v)^2(t,X _{ [0,t]})\,dt\\
	&= h^4 \Big(E_Q \int_0 ^{ T} (\LLf^Q_tu)^2(X)\,dt+o_{h\to 0^+}(h)\Big)\ 
		\Big(E_Q \int_0 ^{ T} (\LLf^Q_tv)^2(X)\,dt+o_{h\to 0^+}(h)\Big),
\end{align*}
leading to
\begin{align*}
E_Q \int_0 ^{ T-h} h ^{ -1}|D_t^h|\,dt
	\le h \|\LLf^Qu(X)\| _{ L^2(\Qb )}\,\|\LLf^Qv(X)\| _{ L^2(\Qb )}
		+o_{h\to 0^+}(h ).
\end{align*}
Putting everything together, we obtain
\begin{multline*}
\Limh E_Q\int_0 ^{ T-h} \Big|E_Q \big[h ^{ -1}\{u( X_{t+h})-u( X_t)\}\{v( X_{t+h})-v( X_t)\}\mid X _{ [0,t]}\big]\\
	-h ^{ -1}\int_t ^{ t+h} E_Q[\Gaf^Q_s(u,v)(X _{ [0,s]})\mid X _{ [0,t]}]\,ds\Big|\,dt =0.
\end{multline*}
On the other hand, by Corollary  \ref{res-14} applied with the  convolution
kernel $k^h=\frac 1h \1_{[-h,0]}$  and $ \mathcal{A}_t= \sigma(X _{ [0,t]})$, under the assumptions
$\Gaf^Q (u)(X),$ $ \Gaf^Q (v)(X)\in L^1(\Qb)$, we obtain 
\begin{align*}
\Limh E_Q\int_0 ^{ T-h} \Big|h ^{ -1}\int_t ^{ t+h} E_Q[\Gaf^Q_s(u,v)(X _{ [0,s]})\mid X _{ [0,t]}]\,ds - \Gaf^Q_t(u,v)(X _{ [0,t]})\Big|\,dt =0.
\end{align*}
The conclusion of the proof of (a) follows from these last two limits.

\Boulette{(b)}
Changing a little bit the previous arguments, in particular using the assumed Markov property of $Q$, the convolution kernel $ k ^{ -h}:= \frac{1}{h}\1 _{ [0,h]}$ instead of $k^h$, and applying Corollary \ref{res-14} with  $ \mathcal{A}_t= \sigma(X _{ t})$, we obtain similarly
\begin{multline*}
\Limh E_Q\int_h ^{ T} \Big|E_Q \big[h ^{ -1}\{u( X_{t-h})-u( X _{ t})\}\{v( X_{t-h})-v( X _{ t})\}\mid X _{t-h}\big]\\
	-h ^{ -1}\int _{ t-h} ^{ t} E_Q[\Gaf^Q_s(u,v)[X _s]\mid X _{t-h}]\,ds\Big|\,dt =0.
\end{multline*}
On the other hand, as in the proof Corollary  \ref{res-14}   we obtain 
\begin{align*}
& E_Q\int_h ^{ T} \Big|h ^{ -1}\int _{ t-h}^{ t} E_Q[\Gaf^Q_s(u,v)[X_s]\mid X _{ t-h}]\,ds - \Gaf^Q _{ t-h}(u,v)[X_{ t-h}]\Big|\,dt \\
 &\le E_Q\int_h ^{ T} \Big|h ^{ -1}\int _{ t-h}^{ t} \Gaf^Q_s(u,v)[X _s]\,ds - \Gaf^Q _{ t-h}(u,v)[X_{t-h}]\Big|\,dt \\
& = E_Q\int_h ^{ T} \Big| k ^{ -h}\ast \Gaf^Q_t(u,v)[X_t] - \Gaf^Q _{ t-h}(u,v)[X _{ t-h}]\Big|\,dt\\
& \le  E_Q\int_h ^{ T} \Big| k ^{ -h}\ast \Gaf^Q_t(u,v)[X_t] - \Gaf^Q _{ t}(u,v)[X_t]\Big|\,dt\\
&\hskip 5cm+
 E_Q\int_h ^{ T} \Big|  \Gaf^Q _{ t}(u,v)[X_t] - \Gaf^Q _{ t-h}(u,v)[X _{ t-h}]\Big|\,dt.
\end{align*}
We know by Lemma \ref{res-09} that $\Limh E_Q\int_h ^{ T} \Big| k ^{ -h}\ast \Gaf^Q_t(u,v)[X_t] - \Gaf^Q _{ t}(u,v)[X_t]\Big|\,dt=0$.
{ With the additional hypothesis that $(t,x)\mapsto \Gaf^Q_t(u,v)[x]$ is continuous, and because $\Gaf^Q(u,v)$ is integrable, we see that
$	
\Limh E_Q\int_h^T\Big|  \Gaf^Q _{ t-h}(u,v)[X _{ t-h}]- \Gaf^Q _{ t}(u,v)[X_t]\Big|\,dt=0.
$	
}
 Putting everything together we arrive at \eqref{eq-44}.
\end{proof}


\begin{corollary}\label{res-25}
Let $Q$ be any  path measure and take any $u,v$ in the class  $\UQ_2$.\\
Then, 
for almost all $t$,
 \begin{align*}
\Limh E_Q \big[h ^{ -1}\{u( X_{t+h})-u( X _{ t})\}\{v( X _{ t+h})-v( X _{ t})\}\big]
	= E_Q \Gaf^Q _{ t}(u,v)(X).
\end{align*} 
If in addition, 
 $\UU\subset \dom\LLb^Q$,
$ \LLb^Qu(X),\LLb^Qv(X)\in L^2(\Qb)$,
$(u,u), (v,v)\in\dom\Gab^Q$ and $\Gab^Q (u)(X),\Gab^Q (v)(X)\in L^1(\Qb)$, 
and the hypotheses of Lemma \ref{res-01}-(b)  are satisfied, 
then  for almost all $t$, 
 \begin{equation}\label{eq-45}
 \begin{split}
 \Limh E_Q \big[h ^{ -1}\{u( X_{t-h})-u( X _{ t})\}&\{v( X _{ t-h})-v( X _{ t})\}\big]
	\\&= E_Q \Gaf^Q _{ t}(u,v)(X)= E_Q \Gab^Q _{ t}(u,v)(X).
 \end{split}
 \end{equation} 
\end{corollary}

\begin{proof}
The first statement follows directly from statement (a) of Lemma  \ref{res-01} with  Fubini and  Jensen.
Our additional hypotheses on $\UU,$ $u$ and $v$ mean that $u$ and $v$ belong to $ \mathcal{U}_2^{ Q^*}.$   Applying (a) to $Q^*$ instead of $Q$, we have
\begin{multline*}
\Limh E_Q\int_h ^{ T} \Big|E_Q \big[h ^{ -1}\{u( X_{t-h})-u( X_t)\}\{v( X_{t-h})-v( X_t)\}\mid X _{ [t,T]}\big]\\
	-   \Gab^Q_t(u,v)(X)\Big|\,dt =0.
\end{multline*}
With Fubini and Jensen again, we see that    (b)  of Lemma \ref{res-01}, and this identity imply  \eqref{eq-45}.
\end{proof}

\subsection*{Integration by parts formula}

The following easy result is pointed out because it is a technical argument of the proof of next Theorem \ref{res-02}.
\begin{lemma}\label{res-15}
For any measurable bounded function $u\in \dom\LLf^Q$ such that $\LLf ^Qu[\Xb] \in L^1(\Qb),$ and all $0\le s\le t\le T,$
\begin{align*}
E_Q[u(\Xb_t)-u(\Xb_s)\mid X_s]
	=E_Q \left[\int_s^t \LLf^Qu(\Xb_r)\,dr\mid X_s\right].
\end{align*}
For any measurable bounded function  $u\in\dom\LLb^Q$ such that $ \LLb^Qu[\Xb]\in L^1(\Qb),$ and all $0\le s\le t\le T,$
\begin{align*}
E_Q[u(\Xb_t)-u(\Xb_s)\mid X_t]	=-E_Q \left[\int_s^t \LLb^Qu(\Xb_r)\,dr\mid X_t\right] .
\end{align*}
\end{lemma}
\begin{proof}
The first equality is obvious. Let us look at the second one:
\begin{align*}
E_Q[&u(\Xb_t)-u(\Xb_s)\mid X_t]
	=E _{ Q^*}[u^*(T-t,X _{ T-t})-u^*(T-s, X _{ T-s})\mid X _{ T-t}]\\
	&=-E _{ Q^*} \left[\int _{T-t} ^{ T-s}\LLf ^{ Q^*}u^*(r,X_r)\, dr\mid X _{ T-t}\right] 
	= -E _{ Q^*} \left[\int _{s} ^{ t}\LLf ^{ Q^*}u^*(T-r,X _{ T-r})\, dr\mid X _{ T-t}\right] \\
	&=- E _{ Q^*} \left[\int _{s} ^{ t}\LLb ^{ Q}u(r,X _{ T-r})\, dr\mid X _{ T-t}\right] 
	= -E _{ Q} \left[\int _{s} ^{ t}\LLb ^{ Q}u(r,X _{ r})\, dr\mid X _{ t}\right] ,
\end{align*}
as announced.
\end{proof}

Next Theorem \ref{res-02} is the cornerstone of the proofs of  time reversal formulas. 
Before stating it, let us introduce some notation.
For any path measure $Q$, we define
\begin{align*}
&\LLf ^Q_t u[X_t]:=E_Q \Big(\LLf^Q_t u( X _{ [0,t]})\mid X_t\Big),
\qquad
\Gaf^Q_t(u,v)[X_t]:=E_Q\Big(\Gaf^Q_t(u,v)(X _{ [0,t]})\mid X_t\Big),
\\
&\LLb^Q_t u[X_t]:=E_Q \Big(\LLb^Q_t u( X _{ [t,T]})\mid X_t\Big),
\qquad
\Gab^Q_t(u,v)[X_t]:=E_Q\Big(\Gab^Q_t(u,v)(X _{ [t,T]})\mid X_t\Big),
\end{align*}
where we use square brackets $[X_t]$ to specify  the conditional expectation knowing $ X_t$, provided it is well defined. Of course, if $Q$ is Markov, then $\LLf ^Q_t u[X_t]=\LLf ^Q_t u(X_t),$ and so on.
We introduce the class of functions
\begin{equation}\label{eq-42}
\UQ := \Big\{ u\in \UU; \LLf^Qu[\sbt]\in L^1(\qb),\ \Gaf^Q (u)[\sbt]\in L^1(\qb)  \Big\}.
\end{equation}
Comparing with \eqref{eq-42b}, we see that the  differences with $\UQ_2$ are  the conditional expectations with respect to $X_t$ and that  $\LLf^Qu[\sbt]$  stands in $L^1(\qb)$ instead of $L^2(\qb)$. The integrability improvement is  useful at Section \ref{sec-rw} and in the companion paper \cite{CL20}  when establishing time reversal formulas for jump processes under a finite entropy hypothesis.

\begin{theorem}[IbP of the carré du champ]\label{res-02}
Let $P\in\MO$ be any path measure. Take two functions $u,v$ in  $\UP.$

\begin{enumerate}[(a)]
\item
If
\begin{align}\label{eq-39a}
u\in\dom \LLb^P\quad \textrm{and}\quad \LLb^Pu(X)\in L^1(\Pb),
 \end{align} 
then for almost every $t$ 
\begin{align}\label{eq-24}
E_P\Big((\LLf_t^Pu+\LLb_t^Pu)[X_t]v(X_t) + \Gaf_t^P(u,v)[X_t] \Big)
=0.
\end{align}

\item
Suppose that  $P$ is Markov, 
\begin{equation}\label{eq-39b}
(t,x)\mapsto \Gaf^P_t(u,v)(x) \textrm{ is continuous},
\end{equation}
the class of functions $\UP$ determines the  weak convergence of Borel measures on $\XX$,
and the linear form 
\begin{equation}\label{eq-39c}
w\in \UPb \mapsto E_P\Iii \Gaf^P_t(u,w_t)(X_t)\,dt
\end{equation}
on $\UPb:= \left\{ w\in C_b(\iX); \ w(t,\sbt)\in \UP,\ \forall 0\le t\le T\right\} $ defines a finite measure on $\iX$.\\
Then, \eqref{eq-39a} holds   and therefore \eqref{eq-24} is satisfied.
\end{enumerate}
\end{theorem}

\begin{remarks}\label{rem-05}\ 
\begin{enumerate}[(a)]
\item
The assumption \eqref{eq-39c} is an integration by parts formula.
\item
Statement (a) is really significant when $P$ is a Markov measure because in this case $\LLf^P[X]=\LLf^P(X)$, $\LLb^P[X]=\LLb^P(X)$ and $\Gaf^P[X]=\Gaf^P(X)$: we do not loose any information and this carries all the necessary material to derive a time reversal formula. We state it in the general form to stress that the Markov property does not play any role in the proof of statement (a).
\item
Using the notion of osmotic extended generator 
\begin{align*}
\LLo P_tu(x):= (\LLf^P_t+\LLb^P_t)u[x]/2,
\end{align*} 
the IbP formula writes as 
\begin{align*}
\IX v\LLo P_t u\, dP_t= -\ud \IX \Gaf ^P_t(u,v)[x]\, P_t(dx)
	= -\ud \IX \Gab ^P_t(u,v)[x]\, P_t(dx),
\end{align*}
where last equality is Corollary \ref{res-25}, provided that the extra hypotheses of this corollary are satisfied.
We see that it extends the usual integration by parts formula stated at Proposition \ref{res-23} below, which is only valid for \emph{stationary} Markov measures. 
\item
The symmetry of the carré du champ implies 
\begin{align*}
\IX v\LLo P_t u\, dP_t=\IX u\LLo P_t v\, dP_t.
\end{align*}

\item
By Proposition \ref{res-10} we know that for any $u\in\dom\LLo P:=\dom\LLf^P\cap\dom\LLb^P$ such that $E_P\Iii (|\LLf^P_tu|+|\LLb^P_tu|)(X)\,dt< \infty,$ the limit
\begin{align*}
\LLo P_t u(X_t) =\Limh \frac 1hE_P\Big( \frac{u(X _{ t+h})+u(X _{ t-h})}{2}-u(X_t)\mid X_t \Big)
\end{align*}
takes place in $L^1(\Pb).$
\end{enumerate}
\end{remarks}

\begin{proof}[Proof of Theorem \ref{res-02}]
We start proving the IbP formula \eqref{eq-24} assuming that $u$ and $v$ belong to $\UP_2,$ and  using both hypotheses \eqref{eq-39a} and \eqref{eq-39b}. Once this is done, we extend the result to the case where $u$ and $v$ are in $\UP.$ Finally, we shall see at the end of the proof that it is a simple matter to remove one assumption among \eqref{eq-39a} and \eqref{eq-39b}.

\par\medskip\noindent $\bullet$\ Proof of \eqref{eq-24} under the hypotheses: $u,v\in\UP_2,$ \eqref{eq-39a} and \eqref{eq-39b}.\quad 
It is based on the elementary identity
\begin{equation}\label{eq-46}
 \begin{split}
[(u _{ t+h}-u_t)+&(u _{ t-h}-u_t)]v_t\\
	&= -(u _{ t-h}-u _{ t})(v _{ t-h}-v _{ t})+v _{ t}(u _{ t+h}-u _{ t})
		-v _{ t-h}(u_t-u _{ t-h}),
\end{split}
\end{equation}
which implies 
\begin{align*}
&E_P\Big(\big\{E_P[u (X _{ t+h})-u(X_t)\mid X_t]+E_P[u(X _{ t-h})-u(X_t)\mid X_t]\big\}\ v(X_t)\Big)\\
	&\  = -E_P\Big(\{u(X  _{ t-h})-u(X _{ t})\}\{v(X _{ t-h} )-v(X_{ t})\}\Big)\\
	&\hskip 0,6cm  +E_P\Big(v(X _{ t})E_P[u(X _{ t+h})-u(X_t)\mid X _{ t}]\Big)
	-E_P\Big(v(X _{ t-h})E_P[u(X_t)-u(X _{ t-h})\mid X _{ t-h}]\Big).
\end{align*}
Dividing both sides by $h>0$, letting $h\to 0^+$, 
\begin{equation}\label{eq-04}
\begin{split}
E_P&[(\LLf_t u+\LLb_t u)[X_t]\, v(X_t)]\\
	&=-E_P \Gaf^P_t(u,v)[X_t]
	 +\Limh E_P\Big(v( X _{ t})E_P[u( X _{ t+h})-u( X_t)\mid  X _{ t}]\Big)\\
	&\hskip 4.5 cm -\Limh E_P\Big(v( X _{ t-h})E_P[u( X_t)-u( X _{ t-h})\mid X _{ t-h}]\Big),
\end{split}
\end{equation}
and the proof will be complete once we show that the last two terms cancel each other.
\\
Let us present  some justifications for \eqref{eq-04}. 
We denote     for any $0<h\le T,$ 
\begin{align*}
\Xb ^h_t&:=\Xb _{ t+h}=(t+h, X _{ t+h}),\quad 0\le t\le T-h,\\
\Xb  ^{ -h}_t&:=\Xb _{ t-h}=(t-h,X _{ t-h}),\quad h\le t\le T.
\end{align*}
Because $u$ is taken in $\UP_2,$ by the first part of  Proposition \ref{res-10} (and Jensen's inequality),   we have $\Limh h ^{ -1}E_P[u(\Xb^h)-u(\Xb)\mid\Xb] =\LLf^Pu[\Xb]$ in $L^2(\Pb ),$ and a fortiori in $L^1(\Pb).$ But $v(\Xb)$ is  a bounded function. Hence, 
\begin{align*}
\Limh E_P\Big(v(\Xb)\, h ^{ -1} [u(\Xb^h)-u(\Xb)]\mid\Xb\Big) =v(\Xb)\LLf^Pu[\Xb]
\textrm{\quad  in }L^1(\Pb ).
\end{align*}
Similarly, with the second part of Proposition \ref{res-10}, under the assumption \eqref{eq-39a}
\begin{align}\label{eq-40a}
\Limh E_P\Big( v(\Xb)\, h ^{ -1} [u(\Xb ^{ -h})-u(\Xb)]\mid\Xb\Big)  =v(\Xb)\LLb^Pu[\Xb]
\textrm{\quad  in }L^1(\Pb ).
\end{align}
With Fubini and Jensen, this proves 
\begin{align}\label{eq-23}
\begin{split}
\Limh \int_h ^{ T-h} \Big|h ^{ -1} E_P\Big([\{u (X _{ t+h})-&u(X_t) \}+\{u(X _{ t-h})-u(X_t)\}]\,v(X_t)\Big)\\&-
E_P \Big((\LLf_t^Pu+\LLb_t^Pu)[X_t]\,v(X_t)\Big) \Big|\, dt=0 .
\end{split}
\end{align}
Under the assumption \eqref{eq-39b} and because  $u$ and $v$ are assumed to belong to $\UP_2$, by Lemma \ref{res-01}-(b) we  have also
\begin{equation} \label{eq-40b}
\begin{split}
\Limh \int_h ^{ T} \Big|E_P \Big(h ^{ -1}\{u( X_{t})-u( X _{ t-h})\}\{v( X _{ t})-v&( X _{ t-h})\}\Big)\\
	&-   E_P \Gaf^P _{ t}(u,v)[X_t]\Big|\,dt =0.
\end{split}
\end{equation} 
It remains to prove that the last terms of \eqref{eq-04} cancel each other by showing that 
\begin{align}
\Limh \int _0^{ T-h}\Big| E_P\Big(v( X _{ t})h ^{ -1}[u( X _{ t+h})-u( X_t)]\Big)
	&-E_P( v(\Xb_t)\LLf u[\Xb_t])\Big|\,dt=0 ,\label{eq-05}\\
\Limh \int_h^T \Big| E_P\Big(v( X _{ t-h})h ^{ -1}[u( X_t)-u( X _{ t-h})]\Big)
	&-E_P\big(v(\Xb_t)\LLf u[\Xb_t])\Big|\,dt=0 .\label{eq-18}
\end{align}
The leftmost integrand  of \eqref{eq-05} is 
$E_P(v(X_t) \{k ^{ h}\ast \LLf u[\Xb]\}_t)$, so that the identity follows because $v(\Xb_t)$ is in $L^\infty(P)$ and $\Limh \{k ^{ h}\ast \LLf u[\Xb]\}=\LLf u[\Xb]$ in $L^1(\Pb)$  by Lemma \ref{res-09}.
\\
On the other hand, \eqref{eq-18} is true because 
\begin{enumerate}[(i)]
\item
$ E_P\Big(v(X _{ t-h})h ^{ -1}[u(X_t)-u(X _{ t-h}) ]\Big)
	= E_P\Big(v(X _{ t-h})\{k ^{ -h}\ast \LLf u[\Xb]\} _{ t}\Big)$;
\item
$\Limh k ^{ -h}\ast\LLf u[\Xb]=\LLf u[\Xb]$ in $L^1(\Pb)$;
	\item
$\Limh v(X _{ t-h})=v(X _{ t^-})=v(X_t), \Pb\ae$ 
\end{enumerate}
Item (i) follows from Lemma \ref{res-15} and \eqref{eq-22},
 (ii) is a direct consequence of Lemma \ref{res-09}, and 
 (iii) follows because the sample paths are left-limited,  it is assumed that     $v$ is continuous and bounded, and $X_t=X _{ t^-}$ for almost every $t$, $P\ae$ because the sample paths are \cadlag.
 \\
We have proved
\begin{align}\label{eq-24b}
E _{ P}\Iii \Big| (\LLf_t ^{ P}u+\LLb_t ^{ P}u)[X_t]v(X_t) + \Gaf_t ^{ P}(u,v)[X_t]  \Big|\,dt=0,
\end{align}
and therefore  \eqref{eq-24}, under the hypotheses: $u,v\in\UP_2,$  \eqref{eq-39a} and \eqref{eq-39b}. Let us relax  this hypothesis by considering functions $u$ and $v$ in $\UP$ instead of $\UP_2.$

\par\medskip\noindent $\bullet$\ Proof of \eqref{eq-24} under the hypotheses: $u,v\in\UP,$ \eqref{eq-39a} and \eqref{eq-39b}.\quad The proof of this extension relies on a localization argument.
For any $u,v\in\UP$ and any $k\ge 1,$ we define the stopping time $$ \tau^k:=\inf\Big\{t\in\ii; \int_0^t|\LLf^P_su(X _{ [0,s]})|\,ds+\int_0^t|\LLf^P_sv(X _{ [0,s]})|\,ds\ge k\Big\}$$
and consider the sequence of stopped path measures $P ^{ k}:=(X ^{ \tau^k})\pf P,$ $k\ge 1.$ Clearly $$\Lim k \tau^k= \infty,\ P\ae$$ because $u$ and $v$ belong to $\dom\LLf^P.$
For any $k\ge 1,$ the functions $u$ and $v$ are in $ \mathcal{U} ^{ P^k}_2,$ therefore we have just proved that $P^k$ verifies  \eqref{eq-24b}: 
\begin{align*}
E _{ P^k}\Iii \Big| (\LLf_t ^{ P^k}u+\LLb_t ^{ P^k}u)[X_t]v(X_t) + \Gaf_t ^{ P^k}(u,v)[X_t]  \Big|\,dt=0.
\end{align*}
On the other hand,
 $\LLf ^{ P^k}_tu[x]= E_P[\1 _{ \{t<\tau^k\}}\LLf^P_tu(X)\mid X_t=x],$ { $\LLb ^{ P^k}_tu[x]=E_P \1 _{ \{t\le \tau^k\}}\LLb^P_tu,$} and  $\Gaf ^{ P^k}_t(u,v)=\1 _{ \{t<\tau^k\}}\Gaf ^{ P}_t(u,v).$ Hence
\begin{align*}
0&=E _{ P^k}\Iii \Big| (\1 _{ \{t<\tau^k\}}\LLf_t ^{ P}u+\1 _{ \{t\le \tau^k\}}\LLb_t ^{ P}u)[X_t]v(X_t) + \1 _{ \{t<\tau^k\}}\Gaf_t ^{ P}(u,v)[X_t] \Big|\,dt\\
&=E _{ P}\Iii \Big| (\1 _{ \{t<\tau^k\}}\LLf_t ^{ P}u+\1 _{ \{t\le \tau^k\}}\LLb_t ^{ P}u)[X_t]v(X_t) + \1 _{ \{t<\tau^k\}}\Gaf_t ^{ P}(u,v)[X_t] \Big|\,dt\\
&=E _{ P}\Iii \Big| (\LLf_t ^{ P}u+\LLb_t ^{ P}u)[X_t]v(X_t) +\Gaf_t ^{ P}(u,v)[X_t] \Big|\,dt.
\end{align*}
The second equality holds because $P$ and $P^k$ match on $ \left\{t\le \tau^k\right\} ,$ and last equality follows  letting $k$ tend to infinity  by dominated convergence under our integrability assumptions.  
We have proved \eqref{eq-24} under the hypotheses: $u,v\in\UP,$  \eqref{eq-39a} and \eqref{eq-39b}.

This proof was based on the convergence of the identity \eqref{eq-46} as $h$ tends to zero. But for this convergence to hold, it is sufficient that only three of its four terms converge. We take advantage of this remark to complete the proof.

\Boulette{(a)}
Let us remove \eqref{eq-39b}. This assumption was used to obtain \eqref{eq-40b} and was not used anywhere else. Hence, the  limits of the other three terms of \eqref{eq-46} are valid even in absence of \eqref{eq-39b}, showing in return that in addition to \eqref{eq-24}, \eqref{eq-40b} holds true.

\Boulette{(b)}
Let us remove \eqref{eq-39a}. This assumption was used to obtain \eqref{eq-40a} and was not used anywhere else. Hence, the  limits of the other three terms of \eqref{eq-46} are valid even in absence of \eqref{eq-39a}, showing in return that in addition to \eqref{eq-24}, the expectation of \eqref{eq-40a} holds true:  the limit 
\begin{align}\label{eq-47}
\Limh E_P\int _{ [h,T]} w _t(X_t) h ^{ -1} \left\{ u(X _{ t-h})-u(X_t)\right\} \,dt=:\overleftarrow{\ell} ^{ P}_u(w )
\end{align}
exists for all $w \in\UPb $ (passing from $\UP$ to $\UPb$ is obvious), and we have
\begin{equation}\label{eq-48}
\overleftarrow{\ell} ^{ P}_u(w)=E_P\Iii \big(-\LLf_t^Pu[X_t]w_t(X_t) - \Gaf_t^P(u,w)[X_t] \big)\,dt,
\qquad w\in \UPb.
\end{equation}
As we  assume that $\UP$ is separating and  $w\mapsto E_P\Iii   \Gaf_t^P(u,w)[X_t] \,dt$ defines a finite measure, $\overleftarrow{\ell} ^{ P}_u$ is also a finite measure on $\iX$ because $\LLf^Pu\in L^1(\Pb).$ It is absolutely continuous with respect to $\pb$ because $\IiX|w|\,d\pb=0$ implies  $\IiX \Gaf_t^P(u,w_t)[x]\,\pb(dtdx)=0.$
Moreover, since $\UP$ is convergence-determining, so is $\UPb$, and with \eqref{eq-47} and Proposition \ref{resh-03b} we see that
$
 \LLb^Pu={d\overleftarrow{\ell}^{ P}_u}/{d\pb}.
$
We conclude with  \eqref{eq-48} that the IbP formula \eqref{eq-24} is satisfied.
 \end{proof}

\begin{remark}
Another very similar proof is based on the elementary identity
\begin{equation*}
\begin{split}
[(u _{ t+h}-u_t)+&(u _{ t-h}-u_t)]v_t\\
	&= -(u _{ t+h}-u _{ t})(v _{ t+h}-v _{ t})+v _{ t}(u _{ t-h}-u _{ t})
		-v _{ t+h}(u_t-u _{ t+h}).
\end{split}
\end{equation*}
Doing this, one sees that (a) is still valid. But (b) does not follow so easily because one cannot drop \eqref{eq-39a}.
\end{remark}

\begin{corollary}\label{res-24}
Let $u$ be in $ \mathcal{U}^P$ and suppose that \eqref{eq-39a} is fulfilled. 
Then, $u\in\dom \LLo P$ and  for almost every $t$
\begin{align*}
\IX \LLo P_t u\,dP_t
=0.
\end{align*}
\end{corollary}

\begin{proof}
Apply Theorem \ref{res-02} with $v=1.$ 
\end{proof}

\subsection*{Stationary Markov measure}\label{sec-statio}
To make the point of Remark \ref{rem-05}-(c) precise,  let us recall what time reversal does with stationary Markov measures.

We consider a \emph{stationary} Markov measure $P\in\MO$ with stationary measure $\mm,$ i.e.\  $P_t=\mm,$ for all $t\in\ii.$ Stationary means that for any real numbers $t_1,\dots,t_k$ and $h$, the laws of $(X _{ t_1+h},\dots,X _{ t_k+h})$ and $(X _{ t_1},\dots,X _{ t_k})$ under $P$ are the same. As $P$ is Markov, it is sufficient that  this property holds for $k= 2.$  We restrict the time interval to $\ii.$

Define the class of functions
\begin{align*}
\mathcal{V}:= \left\{v\in\dom\LLf^P\cap L^2(\mm);\  \LLf^Pv\in L^1(\mm)\cap L^2(\mm)\right\} .
\end{align*}
The  adjoint $(\LLf^P _{ | \mathcal{V}})^*$   in $L^2(\mm)$ of the restriction $\LLf^P _{ | \mathcal{V}}$ to $ \mathcal{V}$ of the forward generator $\LLf^P$  of $P$     is defined by: $\IX v(\LLf^P _{ | \mathcal{V}})^*u\,d\mm=\IX u\LLf^P v\,d\mm,$ for any $u,v\in \mathcal{V}.$

\begin{lemma}\label{res-19}
Suppose  that $P$ is Markov and stationary, then:
$
\LLb^P _{ | \mathcal{V}}=(\LLf^P _{ | \mathcal{V}})^*.
$
\end{lemma}

\begin{proof}
Fix $t,h$ such that $0\le t\le t+h\le T$ and take $u,v\in \mathcal{V}.$ By stationarity
\begin{align*}
E_P[u(X_t)&\{v(X _{ t+h})-v(X_t)\}]\\
	&=  E_P[\{u(X _{ t-h})-u(X _{ t})\}v(X _{ t})]
	+E_P[u(X_t)v(X _{ t+h})-u(X _{ t-h})v(X _{ t})]\\
	&=  E_P[\{u(X _{ t-h})-u(X _{ t})\}v(X _{ t})].
\end{align*}
Dividing by $h>0$ and letting it tend to zero, we conclude with Proposition \ref{res-10}.
\end{proof}
This is a well-known result. One of its versions in the framework of discrete time  was  published by Nelson in 1958  \cite{N58}. We find it pleasant to provide an elementary  proof in the continuous-time setting, based on stochastic derivatives: a tool developed by Nelson himself a decade later.

\begin{proposition}[Integration by parts] \label{res-23}
Suppose  that $P\in\MO$ is Markov and stationary, then for any $u,v\in \mathcal{V},$ such that $uv\in \mathcal{V,}$
\begin{align*}
\IX u  \mathcal{L} ^{ \mathrm{sym},P} v\,d\mm=-\ud \IX \Gaf^P(u,v)\, d\mm,
\end{align*}
where 
\begin{align}\label{eq-Lsym}
\mathcal{L} ^{ \mathrm{sym},P}:=(\LLf^P+(\LLf^P)^*)/2=(\LLf^P+\LLb^P)/2=\LLo P 
\end{align}
is the algebraic symmetrization of $\LLf^P$.
\end{proposition}
 
 \begin{proof}
 Let us denote for simplicity $A:=\LLf^P _{ | \mathcal{V}}.$
 Of course, $\IX Au\,d\mm=0$ because $\IX Au\,d\mm=\IX \1 Au\,d\mm=\IX uA^*\1\,d\mm$ and the stationarity implies that  $A^*\1=0.$ Therefore,
\begin{align*}
\IX \Gaf^P(u,v)\,d\mm
	&=\IX \{A(uv)-uAv-vAu\}\,d\mm\\
	&=-\IX \{uAv+vAu\}\,d\mm
	=-\IX \{uA^*v+vA^*u\}\,d\mm\\
	&=-2\IX u\bar Av\,d\mm
\end{align*}
with $\bar A:=(A+A^*)/2.$ 
We conclude with Lemma \ref{res-19}.
 \end{proof}

\subsection*{Finite entropy}

Up to now the entropy did not play any role. Let us write some words about  it in preparation to forthcoming time reversal formulas. 
\\
Comparing statements (a) and (b) of Theorem \ref{res-02}, we see that (b) is  easier to verify than (a), because (a) requires that $u$ is  in the domain of the backward generator: a property which is not known a priori. On the other hand,  the assumption \eqref{eq-39b} in (b) is too much demanding for some applications we have in mind, where a finite entropy condition destroys this regularity  in presence of jumps, see \cite{CL20}. 

We are going to investigate time reversal of Markov measures $P$ verifying the finite entropy condition \eqref{eq-36}: 
$
H(P|R)< \infty,
$
where the time reversal $R^*$ of a reference Markov measure $R$ is accessible via Theorem \ref{res-02}-(b). Then, taking advantage of the elementary identity $H(P^*|R^*)=H(P|R)< \infty,$ a deep insight of Föllmer already encountered at \eqref{eq-38},  we shall be in position to build a large enough class $\UP$ and to verify the assumptions of Theorem \ref{res-02}-(a) for  $P$.

\section{Time reversal of a diffusion process in $\ZZ$}
\label{sec-diff1}

In this section, the IbP formula of Theorem \ref{res-02} is used to obtain at Theorem \ref{res-08} a  time reversal formula for diffusion measures.

\subsection*{Reference diffusion measure}

The path space is the set $\OO=C([0,T],\ZZ)$ of all continuous trajectories from $[0,T]$ to $\ZZ.$  The main reference measure we have in mind is the reversible Kolmogorov diffusion $R$ defined at \eqref{eq-09}.

\subsection*{Finite entropy in a diffusion setting}
Take $Q\in\PO$ such that
\begin{align}\label{eq-10}
H(Q|R)< \infty.
\end{align}
We know by the Girsanov theory under a finite entropy  condition   \cite{Leo11a}, that when  $R$ fulfils  the uniqueness condition:
\begin{align}\label{eq-H17a}
\forall R'\in\MO,\ [R'\in\MP(R_0,\aa,b^R) \textrm{ and }R'\ll R]\implies R'=R,
\end{align}
 there exists some $\ZZ$-valued predictable   process $ \bbb QR$ which is defined $\Qb\ae$    such that   $Q$ solves the martingale problem
\begin{align}\label{eq-11}
Q\in\MP(Q_0,\vam+\aa  \bbb QR ,\aa ).
\end{align}
Recall Remark \ref{rem-UG}  for a  setting where the uniqueness condition  \eqref{eq-H17a} is satisfied.
\\
Furthermore, because of the uniqueness of the solution to $\MP(\mm,\vam,\aa)$, we know that
\begin{align*}
\frac{dQ}{dR}
	&= 	\1 _{ \left\{ dQ/dR>0\right\} }\ 
		\frac{dQ_0}{dR_0}(X_0)
		\exp \left( \Iii \bbb QR_t\scal dM^R_t
		-\Iii |\bbb  Q R  _t|^2 _{ \aa(X_t)} /2\ dt\right)\nonumber\\
	&= 	\1 _{ \left\{ dQ/dR>0\right\} }\ 
		\frac{d Q_0}{dR_0}(X_0)
		\exp \left( \Iii \bbb QR_t\scal dM^Q_t
		+\Iii |\bbb  Q R  _t|^2 _{ \aa(X_t)} /2\ dt\right) ,
\end{align*}
where 
\[
dM^R_t=dX_t-\vam(X_t)\,dt
\quad \textrm{and}\quad dM^Q=dX_t-(\vam(X_t)+\aa(X_t) \bbb QR_t)\,dt,
\]
and we denote $$| \beta|^2_\aa:= \beta\scal \aa \beta.$$
Moreover, 
\begin{align}\label{eq-12}
H( Q|R)
	= H( Q_0|R_0)+E_ Q  \Iii  |\bbb  Q R  _t|^2 _{ \aa(X_t)}/2\ dt.
\end{align}
Of course, in view of this identity,  $H( Q |R )< \infty$ implies that $E_ Q  \Iii   |\bbb  Q R  _t|^2 _{ \aa(X_t)}\, dt$ is finite. 

\begin{claim}
If in addition $ Q $ is Markov, then the process $ \bbb  Q R  $ turns out to be  a vector field: 
\begin{align*}
\bbb  Q R  _t=\bbb  Q R  (\Xb_t), \ \Qb  \ae
\end{align*} 
\end{claim}

\begin{proof}
Indeed, we see with \eqref{eq-11} that
\begin{align*}
[\vam(X_t)+\aa(X_t) \bbb  Q R  _t]\,dt 
	&=E_ Q ( dX_t\mid X _{ [0,t]})
	=E_ Q ( dX_t\mid X _t)\\
	&=\vam(X_t)+\aa(X_t) E_ Q (\bbb  Q R  _t\mid X_t)]\,dt ,
\quad  Q \ae
\end{align*}
Remark that all the above conditional expectations are well-defined; in particular $E_ Q (\bbb  Q R  _t\mid X_t)$ is meaningful because of \eqref{eq-12} and the finite entropy assumption \eqref{eq-10}.
It follows that for all $0\le t\le T,$
	$\aa(X_t) \bbb  Q R  _t=\aa(X_t) E_ Q (\bbb  Q R  _t\mid X_t),$ $\Qb  \ae$ 
\end{proof}

Moreover, we observe that 
\begin{align*}
H(Q|R)-H(Q_0|R_0)
	=H(Q|R ^{ Q_0})
	=E_Q\Iii\ud |\vf ^{ Q|R}|^2_\gg (\Xb_t)\, dt
\end{align*}
is an average kinetic action, 
where $\vf ^{ Q|R}:=\aa \bbb QR$ and 
\begin{align*}
\gg=\aa ^{ -1}. 
\end{align*}


\subsection*{Nelson's velocities}

The \emph{forward stochastic velocity} $\vf^Q$ is 
\begin{align*}
\vf^Q(t,x):=\Lf^Q_t[\Id](x)=  \Limh E_Q \Big( \frac{X _{ t+h}-X_t}{h}\mid X_t=x\Big),
\end{align*}
and similarly, we define the \emph{backward velocity}
\begin{align*}
\vb^Q(t,x):=\Lb^Q_t[\Id](x)= \Limh E_Q \Big( \frac{X _{ t-h}-X_t}{h}\mid X_t=x\Big),
\end{align*}
whenever these expressions are meaningful. These velocities might not be well defined because of a lack of integrability. However, under a finite entropy condition,  Proposition \ref{res-05} below tells us that  they are well defined in the setting we are interested in.
\\
The \emph{current velocity} is
\begin{align*}
\vc Q:=  (\vf^Q-\vb^Q)/2
\end{align*}and the \emph{osmotic velocity} is
\begin{align*}
\vo Q:=(\vf^Q+\vb^Q)/2.
\end{align*}
We immediately observe that
\begin{align*}
\left\{ \begin{array}{lcl}
\vf&=&\ \ \, \vcu+\vos,\\ \vb&=&-\vcu+\vos
\end{array}\right.
\qquad \textrm{and}\qquad
\left\{ \begin{array}{lcl}
\vc{Q^*}_t&=&-\vc Q _{ T-t},
\\ 
\vo {Q^*}_t&=&\ \ \,\vo Q _{ T-t}.
\end{array}\right.
\end{align*}

\subsection*{Entropy under time reversal}

Next result is a central observation in Föllmer's approach to time reversal.

\begin{proposition}\label{res-05}
Under the Hypotheses \ref{ass-01}, let $P$ be a Markov probability measure  such that $H(P|R)< \infty.$\\
  Then, 
there exist two measurable vector fields $\bf ^{ P|R}$ and $\bb ^{ P|R}$ such that 
\begin{align*}
\begin{array}{rcrl}
\LLf^P&=& \partial_t+\vf^P\scal\nabla+ \Delta_\aa /2, &\textrm{where}\quad \vf^P=\vam+\aa \bf ^{ P|R}\\
\LLb^P&=& -\partial_t+\vb^P\scal\nabla+ \Delta_\aa /2,
& \textrm{where}\quad \vb^P=\vam+\aa \bb ^{ P|R}
\end{array}
\end{align*}
with
\begin{align*}
E_P\Iii (|\bf ^{ P|R}|_\aa ^2+|\bb ^{ P|R}|^2_\aa )(\Xb_t)\, dt< \infty,
\end{align*}
and
\begin{align*}
\begin{split}
H(P|R)=H(P_0|R_0)+E _{ P} \Iii \ud & |\bf ^{ P|R}|^2_\aa(\Xb_t)\,dt\\
	&= E _{ P} \Iii \ud |\bb ^{ P|R}|^2_\aa(\Xb_t)\,dt+ H(P_T|R_T).
\end{split}
\end{align*}
\end{proposition}

\begin{proof}
Since $P$ is Markov, so is $P^*:=(X^*)\pf P$.
As the  time reversal mapping $X^*$ is one-one, we have  $H(P|R)=H(P^*|R^*)$.  Hence,
\begin{align*}
H(P|R)=H(P^*|R^*)=H(P^*|R)< \infty,
\end{align*}
where last equality comes from the reversibility of $R$ which implies $R^*=R.$
Again, by Girsanov theory we know that  there is some previsible vector field $ \bf ^{ P^*|R} $ such that
  $P^*$ solves the martingale problem 
$	
\textrm{MP}(P_T,\vam+\aa\bf ^{ P^*|R},\aa).
$	
Denoting
$	
\bb^{P|R}(t,z):= \bf ^{ P^*|R}(T-t,z),
$	
we see that
\begin{align*}
H(P^*|R)&=H(P^*_0|R_0)+E _{ P^*} \Iii \ud |\bf ^{ P^*|R}|^2 _{ \aa}(\Xb_t)\,dt\\
	&= H(P_T|R_T)+ E _{ P} \Iii \ud |\bb ^{ P|R}|^2_\aa(\Xb_t)\,dt,
\end{align*}
as announced. 
\end{proof}

\subsection*{Continuity equation}

Proposition \ref{res-07} below, which is the object of this subsection is not directly linked to time reversal (it is rather complementary). Nevertheless, we present its easy short proof because all the preliminary notions which are needed to its statement and proof appear in the last previous pages.
  
\begin{lemma}\label{res-06}
Under the Hypotheses \ref{ass-01}, let $P\in\PO$ be Markov and such that $H(P|R)< \infty.$ 
Then, any compactly supported function $u\in C ^{ 1,2}_c(\iZ)$ stands in the domain of both $\LLf^P$ and $\LLb^P,$ and $E _{ \Pb} |\LLf^P u(\Xb)|^2< \infty$, $E _{ \Pb} |\LLb^P u(\Xb)|^2< \infty.$
 \\
Moreover  
$
u(\Xb _t)-u(\Xb_0)-\int_0^t \LLf^P _su[X_s]\,ds
$
and 
$
u(\Xb _t)-u(\Xb_T)-\int_t^T \LLb^P _su[X_s]\,ds
$
are respectively genuine (rather than local) forward and backward $P$-martingales. 
\end{lemma}

\begin{proof}
The proofs of the statements concerning $\LLf^P$ and $\LLb^P$ being similar, we focus on $\LLf^P$.
Take $u$ in $C ^{ 1,2}_c(\iZ)$. 
All we have to show is
\begin{align*}
E _{ \Pb} |\LLf^P u(\Xb)|^2< \infty.
\end{align*}
By Proposition \ref{res-05}, 
$$
\LLf ^Pu
	= \partial_t u+ \vam\scal \nabla u+ \ud \Delta_\aa  u+\aa  \bf ^{ P|R}\scal \nabla u
$$ 
and $E _{ \Pb }\big(|\bf ^{ P|R}|_\aa ^2(\Xb)\big)< \infty$. 
Since  $\aa$ and $\vam$ are   locally bounded, 
$\LLf^Ru=\partial_t u+\vam\scal \nabla u+ \ud \Delta_\aa  u$ is bounded.
 The last term is controlled by
\begin{align*}
E _{ \Pb }\big(|\aa \bf ^{ P|R}\scal \nabla u(\Xb)|^2\big)
	\le E _{ \Pb }\big(|\bf ^{ P|R}|_\aa ^2(\Xb)\, |\nabla u|_\aa ^2(\Xb)\big)
	\le \sup |\nabla u|_\aa ^2 { E _{ \Pb }\big(|\bf ^{ P|R}|_\aa ^2(\Xb)\big)}
	< \infty.
\end{align*}
\end{proof}

For any measure $m$ and vector field $ \mathsf{w}$ on $\ZZ,$ we  define $\dive_m(  \mathsf{w})$ by: $$\IZ u \dive_m( \mathsf{w}) \, dm
	:=-\IZ \nabla u\scal \mathsf{w} \,dm,\quad u\in C^1_c(\ZZ),$$ whenever the second integral is meaningful.

As a consequence of Proposition \ref{res-05}, we obtain

\begin{proposition}[Continuity equation] \label{res-07}
Under the Hypotheses \ref{ass-01}, let $P\in\PO$ be Markov and such that $H(P|R)< \infty.$ 
 Then for any $t$, $P_t\ll \mm $  and 
\begin{align*}
\rho_t:= \frac{dP_t}{d\mm }
\end{align*}
solves, in the sense of distributions, the continuity equation
\begin{align*}
\partial_t \rho+\dive _{ \mm }( \rho\vc P)=0.
\end{align*}
Similarly the density 
\begin{align*}
\mu_t:= \frac{dP_t}{d\Leb}
\end{align*}
solves, in the sense of distributions, the continuity equation
\begin{align*}
\partial_t \mu+\dive( \mu\vc P)=0.
\end{align*}
\end{proposition}

\begin{proof}
By Lemma \ref{res-06}, for any $0\le s\le t,$ and any $u\in C ^{ 1,2}_c((0,T)\times \ZZ),$ we have $$E_P[u(\Xb_t)-u(\Xb_s)]
	=\int_s^t E_P \LLf^Pu(\Xb_r)\,dr,$$ 
and similarly, with the definition \eqref{eq-F05} of $\LLb^P$ 
\begin{multline*}
E_P[u(\Xb_t)-u(\Xb_s)]
	=E_{P^*}[u^*(\Xb _{ T-t})-u^*(\Xb _{ T-s})] \\
	=E _{ P^*} \int _{ T-s} ^{ T-t} \LLf ^{ P^*}u^*(\Xb_r)\,dr
	=E_P\int _{ s} ^{ t}  \LLb^Pu(\Xb_r)\,dr.
\end{multline*}
With the expressions of $\LLf^P$ and $\LLb^P$ stated at Proposition \ref{res-05}, this leads us to
\begin{multline*}
0=E_P\int_s^t \ud [\LLf^P-\LLb^P]u(\Xb_r)\,dr
	=E_P\int_s^t [ \partial_r+\vc P\scal\nabla]u(\Xb_r)\,dr\\
	=  \int _{ [s,t]\times\ZZ} [ \partial_r u+\vc P\scal\nabla u] (r,x)\,\rho_r(x)\ \mm(dx)dr
\end{multline*}
which is the first announced continuity equation. The second one follows replacing $ \rho_r(x)\,\mm(dx)$ by $ \mu_r(x)dx.$
\end{proof}

\subsection*{Time reversal formula}

The main result of this section is the following

\begin{theorem}[Time reversal formula]\label{res-08}
Under the Hypotheses \ref{ass-01} on $R$ given at \eqref{eq-09}, let $P\in\PO$ be Markov and such that $H(P|R)< \infty.$ 
\\
Then, the time reversal $P^*$ of $P$ is a solution of the martingale problem 
\begin{align*}
P^*\in \MP(\vf^{P^*},\aa)
\end{align*} 
with 
\begin{align}\label{eq-13d}
\vf^{P^*}_t=\vb^{P}_{T-t}(x)	
	=-\vf^P_{T-t}(x)+   \nabla\scal( \mu_{T-t} \aa)(x)/ \mu_{T-t}(x),
	\quad dtP_t(dx)\ae
\end{align}
where the divergence is in the sense of distributions, $ \mu_t:=dP_t/d\Leb$ and $\vb^P_t$ is defined at almost all $t$.
\\
Furthermore, $P^*$ is the unique solution of $\MP(\vf^{P^*},\aa)$ among the set of  all $Q\in\PO$ such that $H(Q|R) < \infty.$
\\
Denoting $ \rho_t:=dP_t/d\mm $ and $\bo{P|R}:= \bo P-\bo R$, \eqref{eq-13d} is equivalent to   
\begin{align}\label{eq-14c}
\bo{P|R}_t(x)=\nabla \log\sqrt{ \rho_t}(x),
\qquad dtP_t(dx)\ae
\end{align}
where the derivative is distributional and
\begin{align}\label{eq-15b}
\IiZ |\nabla \log \rho_t|^2_\aa \, dP_tdt< \infty.
\end{align}
\end{theorem}

\begin{remarks}[about Theorem \ref{res-08}]\ \begin{enumerate}[(a)] \label{rem-01}
\item
As $H(P|R)< \infty,$ $P_t\ll\mm\ll\Leb$ for all $t$. Hence $ \mu$ and $\rho$ are well defined.
\item
With $\vo{P|R}:= \vo P-\vo R$, this immediately implies that $ P_t\ae$, for almost all $t$,
\begin{align}\label{eq-14b}
\vo{P|R}_t&=\aa  \nabla \log\sqrt{ \rho_t}, \quad
\\\label{eq-13}
 \mu_t\ \vo P_t
 	&= \nabla\scal ( \mu_t \aa )/2,\\
\label{eq-14}
\rho_t \vo{P|R}_t&=\aa  \nabla { \rho_t}/2, 
\end{align} 
in the sense of distributions.

\item
Of course, \eqref{eq-13d} or \eqref{eq-13} are equivalent to
\begin{align}\label{eq-16}
\vb^P_t=
	-\vf^P_t+\nabla\scal \aa +\aa \nabla\log \mu_t,
	\quad P_t\ae
\end{align}
or 
\begin{align}\label{eq-13b}
 \vo P_t=\nabla\scal \aa /2+\aa \nabla \log { \sqrt{\mu_t}},\quad P_t\ae,
\end{align}

\item
The restriction ``$P_t\ae$'' in \eqref{eq-13d}, \eqref{eq-14b}, \eqref{eq-16} and \eqref{eq-13b} prevents $ \rho_t$ and $ \mu_t$ from  vanishing, so that the $\log$ is well defined.

\item
The reference measure 
$$
R^o\in\MP(\Leb, \nabla\scal \aa/2,\aa )
$$
 is the law of a stationary diffusion process with Lebesgue measure as stationary measure ($U=0$). Its forward generator is 
\begin{align*}
\partial_tu+( \nabla\scal \aa)\scal \nabla u/2+ \Delta_\aa u/2
	= \partial_tu+\nabla\scal (\aa\nabla u)/2,\qquad u\in C ^{ 1,2}(\iZ).
\end{align*}
Choosing this reference measure, we see that \eqref{eq-13b} writes as 
$\vo P=\vf ^{ R^o}+\vo {P|R^o}$ with $\vf ^{ R^o}=\nabla\scal \aa/2$, and  $\vo {P|R^o}=\aa \nabla\log \sqrt{ \mu} $ which  is \eqref{eq-14b} with $ \rho_t= \mu_t$ for each $t$ since $\mm ^o=\Leb.$ 

\end{enumerate}\end{remarks}

\begin{proof}[Proof of Theorem \ref{res-08}]
Again, remark that the class of functions $ \mathcal{U}=C ^{ 1,2}_c(\iZ)$ for which  Itô's formula is valid  is an algebra, as required by the hypotheses of  the IbP formula (Theorem \ref{res-02}).
\\
By Proposition \ref{res-05}, we know that there exist two vector fields $\bf^ {P|R }$, $\bb^ {P|R }$ such that  $ \LLf^P= \partial_t+(\vam +\aa \bf^ {P|R })\scal  \nabla +  \Delta_\aa /2 $ and $ \LLb^P= -\partial_t+(\vam +\aa \bb^ {P|R })\scal  \nabla +  \Delta_\aa /2 $ with\begin{align}\label{eq-17}
E_P\Iii |\bo {P|R }|^2_\aa(X_t)\,dt< \infty,
\end{align}
where $$\bos:=(\bf+\bb)/2.$$
 Then, for any test  function $w\in C ^{ 2}_c(\ZZ)$ and almost all $t$, we have
\begin{align}\label{eq-19}
E_P\Big( w(X_t) \vo P_i(\Xb_t) + \aa^i(\Xb_t)\scal \nabla w(X_t)/2 \Big)
=0,
\qquad 1\le i\le n,
\end{align}
where $\vo P_i$ is the $i$-th component of $\vo P,$ and $\aa^i$ is the $i$-th column of $\aa$. This
 follows from an application of Theorem \ref{res-02} with $u(t,x)= \mathrm{proj}_i(x) \,\chi(x),$ $1\le i\le n,$ (where $\mathrm{proj}_i(x):=x_i$ and  $\chi\in C_c ^{ 2}(\ZZ)$ has a compact support  and is equal to $1$ on $\supp w$) and is allowed by 
 Lemma \ref{res-06} which ensures that  $\vf^P_i:= \LLf ^{ P} (\mathrm{proj}_i)= (\vam +\aa \bf ^{P|R })_i$ is in $L^2 _{ \textrm{loc}}(\pb ).$  
Similar estimates hold for $\vb^P.$ 
\\  
Integrating by parts in \eqref{eq-19}, we see that for any compactly supported test function $w$ on $\ZZ$ and almost every $t$,
\begin{align*}
0=\IZ w \, \vo P_i \mu\,d\Leb +\ud\IZ \aa^i\scal \nabla w \,\mu\,d\Leb
	= \IZ w\, \vo P_i \mu\, d\Leb -\ud\IZ w\nabla\scal( \mu \aa^i)\,d\Leb,
\end{align*}
where we drop the time dependence. This proves \eqref{eq-13}.
\\
Let us look at \eqref{eq-14}. Apply \eqref{eq-13} to $P=R$ to obtain
\begin{align}\label{eq-13c}
m\ \vo R
 	= \nabla\scal ( m \aa)/2,
\end{align}
where $ m:=d\mm /d\Leb=e ^{ -U}.$ Although $R$, unlike $P$, might not be a probability measure, it is easy to see that the proof of \eqref{eq-13} directly works with $R$ instead of $P$ (in particular $ \beta ^{ R|R}=0$). 
Because 
$	
\mu= \rho m,
$	
we obtain
$	
\nabla\scal ( \mu \aa)
	=\nabla\scal ( \rho m \aa)
	= m \aa \nabla \rho+ \rho\nabla\scal ( m \aa).
$	
It is important to note that both $\nabla \scal ( \mu \aa)$ and $\nabla \scal ( m \aa)$ are well defined in the sense of distributions (as divergence terms) and are functions by \eqref{eq-13} and   \eqref{eq-13c} (the existence of $\vo P$ is a direct consequence of the assumption that $H(P|R)< \infty$). It follows that 
\[
m \aa \nabla \rho= \nabla\scal( \mu \aa)- \rho\nabla\scal( m \aa)
\] 
is also well defined in the sense of distributions and a function. Putting everything together,
\begin{align*}
 \rho(\vo P-\vo R)
 	= m ^{ -1}(\mu\vo P -  \rho m \vo R)
	=   m ^{ -1}(\nabla\scal ( \mu \aa)-  \rho\, \nabla\scal (  m \aa))/2
	=\aa \nabla \rho/2,
\end{align*}
which is \eqref{eq-14}, and implies \eqref{eq-14c}.
 \\
Finally,  the estimate \eqref{eq-15b} is a rewriting of \eqref{eq-17},  and \eqref{eq-13d} follows directly from \eqref{eq-14c}.
\end{proof}

\section{Time reversal of a diffusion process. Abstract setting} \label{sec-diff2}

We use the IbP formula (Theorem \ref{res-02}) again, to extend at Theorem \ref{res-20} the time reversal formula of Theorem  \ref{res-08} to an abstract diffusion setting where the configuration space $\XX$ is a non-specified Polish space. To our opinion, the main interest of this result is not the extension to an abstract space, but its set of assumptions which sheds light on the close to minimal hypotheses that are necessary for the time reversal formula to hold in a diffusion setting.

\subsection*{Stationary diffusion reference measure}

Assume that the Markov measure $R\in\MO$ is stationary (see page  \pageref{sec-statio}) and in addition  that is a \emph{diffusion} path measure with a  Polish space $\XX$ as its configuration space. In this abstract setting, being a diffusion means that the derivation identity
\begin{align}\label{eq-26}
\Gamma(u,vw)= v\Gamma(u,w)+w \Gamma(u,v)
\end{align}
is valid, and that  for any $P\in\PO$ such that $P\ll R,$ we have 
\begin{align}\label{eq-34}
\Gaf^P=\Gab^P=\Gaf^R= \Gab^R=: \Gamma.
\end{align}
These identities fail in presence of jumps.

\begin{lemma}
Let $R\in\MO$ be an $\mm$-stationary  diffusion path measure with osmotic generator $\LLo R.$ 
For  any  functions $ \rho,u,v\in \mathcal{V}$ such that $ \rho u, uv\in \mathcal{V},$ we have
\begin{align}\label{eq-20} 
\IX \Gamma(\rho,u)v\,d\mm
	=-\IX \left\{ \Gamma(u,v)+2v\LLo Ru\right\} \, \rho\,d\mm. 
\end{align}
\end{lemma}

\begin{proof}
The integration by parts formula is $\IX \Gamma(u,v)\,d\mm=-2\IX v\LLo R u\, d\mm.$
With
$	
\Gamma(u,v)+2v\LLo R u
	= \LLo R (uv)-u\LLo R v-v\LLo R u+2v\LLo R u
	=\LLo R (uv)-u\LLo R v+v\LLo R u,
$	
and the derivation identity \eqref{eq-26},
we obtain
\begin{align*}
-\IX \{ \Gamma(u,v)&+2v\LLo R u\} \, \rho\,d\mm
	=\IX \left\{ -\rho \LLo R (uv)+u \rho \LLo R v-v \rho \LLo R u\right\} \,d\mm\\
	&=\ud \IX \left\{ \Gamma (\rho,uv) - \Gamma(u \rho, v) + \Gamma(v \rho,u)
		\right\} \,d\mm
	=\IX  \Gamma( \rho,u)v  \,d\mm,
\end{align*}
as announced.
\end{proof}

\subsection*{Time reversal formula}

The left hand side of \eqref{eq-20} requires that  for any $u\in \mathcal{V},$ the couple $ (\rho,u)$ stands in the  domain of definition of  $ \Gamma,$ while no regularity of $ \rho$ is needed for having a meaningful right hand side. This suggests the following notion, in the spirit of the definition of a distribution.

\begin{hypothesis}\label{ass-04}
We assume that there exists some algebra $ \mathcal{U}\subset \mathcal{V}$ which is total in $L^2(\mm).$ This means that  for all $u,v\in \mathcal{U}$ we have  $uv\in \mathcal{U}$ and that  for any $w\in L^2(\mm)$,
$\IX u w\,d\mm=0, \ \forall u\in \mathcal{U}$ implies that  $ w=0.$ \\
In addition we suppose that for any $u,v\in \mathcal{V},$ $ \Gamma(u,v)\,\mm$ and $u \LLo R v\,\mm$ are bounded measures.
\end{hypothesis}

Recall Remark \ref{rem-04} for the relevance of this hypothesis.

\begin{definition}\label{def-04}
Let $ \mathcal{U}$ be as in Assumption \ref{ass-04} and let 
  $ \rho$ be a nonnegative measurable function which is defined $\mm\ae$ 
We define the linear operator $ \Gamma( \rho,\sbt)$ on $ \mathcal{U}$ in the weak sense, by the identity \eqref{eq-20}, seeing $(u,v)\mapsto \IX \Gamma( \rho,u)v\,d\mm$ as a bilinear form.  
\end{definition}

\begin{theorem}\label{res-20}
 Let $R\in \MO$ be an $\mm$-stationary diffusion measure: i.e.\  \eqref{eq-26} and \eqref{eq-34} hold, such that the Hypothesis \ref{ass-04} is satisfied, and  for any $u,v\in \mathcal{U},$ $\Gamma(u,v)$ is bounded. Let $P\in\PO$ be  a Markov measure  such that $P\ll R,$   $ \mathcal{U}\subset \dom \LLf^P\cap\dom\LLb^P$, and for any $u\in \mathcal{U}$, $\LLf^P u, \LLb^P u\in L^2(\Pb).$ Then, for any $u\in \UU,$
\begin{align*}
\LLo P u=\LLo R u+ \frac{\Gamma( \rho,u)}{2 \rho}
	=\LLo R u+ \frac{\Gamma( \sqrt\rho,u)}{\sqrt \rho},\qquad dtdP_t\ae,
\end{align*}
where $ \rho_t:=dP_t/d\mm$,  the linear operators $ \Gamma( \rho,\sbt)$ and $ \Gamma( \sqrt\rho,\sbt)$ are defined in the weak sense of Definition \ref{def-04}, and
\begin{align*}
\LLo R=\LL ^{R, \mathrm{sym}}:= (\LLf^R+(\LLf^{R})^\ast)/2,
\end{align*}
 is   the symmetrized extended generator of $\LLf^R,$ see \eqref{eq-Lsym}.
In other words,
\begin{align*}
P^*\in \MP(\LLf^{P^*},\UU)
\end{align*}
where for any $u\in\UU,$
\begin{align*}
\LLf^{P^*}_tu=\LL ^{R, \mathrm{sym}}_{T-t}u+\frac{\Gamma( \sqrt{\rho_{T-t}},u)}{\sqrt {\rho_{T-t}}},\qquad dtdP_t\ae
\end{align*}
\end{theorem}

\begin{proof}
The hypotheses of the IbP formula: Theorem \ref{res-02}, are fulfilled, allowing us to write for all $u,v\in \mathcal{U}$ and almost all $t$ 
\begin{align*}
\IX v_t\LLo P_t u_t\ \rho_t\,d\mm&=-\ud\IZ \Gamma(u_t,v_t)\ \rho_t\,d\mm\\
&=\IX [v_t\LLo R u_t \, \rho_t+ \Gamma( \rho_t,u_t)v_t/2]\,d\mm
\end{align*}
where last equality is   \eqref{eq-20}. The second equality in the first displayed formula follows with \eqref{eq-26} which implies that $ \Gamma( \rho,u)=2 \sqrt \rho \,\Gamma( \sqrt \rho,u).$ The identification of the osmotic and symmetrized generators of a stationary path measure is \eqref{eq-Lsym}. 
\end{proof}

As a direct corollary of this result, we see that any path measure $P\in\PO$ verifying the hypotheses of Theorem \ref{res-20} and such that $P_t=\mm$ for all $0\le t\le T$ (it might not be stationary),   shares its osmotic generator with the $\mm$-stationary path measure $R:$ $\LLo P=\LLo R$ in restriction to $\UU,$ because $ \rho=1.$ 

\section{Current-osmosis decomposition} \label{sec-cod}

This section presents an application of the time reversal formula for a diffusion to entropic optimal transport. The motivation  for a decomposition of the relative entropy into the sum of current and osmotic terms was put forward in the introduction of the paper, see \eqref{eq-35a} and \eqref{eq-35b}. This result is Proposition \ref{res-12} below.

We go back to the setting of Section \ref{sec-diff1} and take the same reference path measure $R$ satisfying the  Hypotheses \ref{ass-01}. 
For any $ \mu_0\in\PZ$ such that $ \mu_0\ll \Leb,$ we denote 
$$
R ^{ \mu_0}(\sbt):=\IZ R(\sbt\mid X_0=x_o)\, \mu_0(dx_o),
$$  
the Markov measure with the same forward dynamics as $R$, i.e.\ $ \LLf ^{ R ^{ \mu_0}}=\LLf^R,$ but with $ \mu_0$ as its  initial marginal.

\begin{definitions}\ 
\begin{enumerate}[(1)]

\item
\emph{(Free energy).}
The \emph{free energy} is defined by
\begin{align*}
\mathcal{F}( \mu):= H(\mu|\mm)/2,
\qquad \mu\in\PZ.
\end{align*}

\item
\emph{(Fisher information).}
It is defined by
\begin{align*}
\mathcal{I}_\aa( \mu|\mm)
	:= \IZ | \nabla \log \sqrt{ d \mu/d\mm}|^2_\aa /2 \ d\mu\in[0, \infty],
\end{align*}
for any $ \mu\in\PZ$ such that $\nabla \log { d \mu/d\mm}$ is well defined in the sense of distributions, and $+ \infty$ otherwise.
\end{enumerate}
\end{definitions}

\begin{proposition}\label{res-12}
Under the hypotheses of Theorem \ref{res-08}, for any $0\le t\le T,$
\begin{align*}
H(P _{ [0,t]}|R ^{ P_0} _{ [0,t]})
	=\mathcal{F}(P_t)- \mathcal{F}(P_0)
	+\int _{ [0,t]}  \Big\{ \left\langle  | \vc {P|R}|^2_\gg /2,P_s \right\rangle 
		+ \mathcal{I}_\aa(P_s|\mm)\Big\}   \,ds.
\end{align*}
\end{proposition}

\begin{proof}
Applying Proposition \ref{res-05}, we see that
\begin{align*}
H(P|R)&= H(P_0|\mm)+E _{ P} \Iii   |\bf ^{ P|R}|^2_\aa(X_t)/2\,dt
	= E _{ P} \Iii  |\bb ^{ P|R}|^2_\aa(X_t)/2\,dt+ H(P_T|\mm)\\
	&= \ud \big(H(P_0|\mm)+ H(P_T|\mm)\big)+\ud E _{ P} \Iii \big(  |\bf ^{ P|R}|^2_\aa +\bb ^{ P|R}|^2_\aa \big)(\Xb_t)/2\,dt.
\end{align*}
On the other hand, with the additive decomposition of the relative entropy
\begin{align*}
H(P|R)=H(P_0|R_0)+\IZ H\big(P(\sbt|X_0=x)\big|R(\sbt|X_0=x)\big)\, P_0(dx),
\end{align*}
we obtain  
\begin{align*}
H(P |R^{ P_0})
	&=H(P_0|P_0)+\IZ H\big(P(\sbt|X_0=x)\big|R(\sbt|X_0=x)\big)\, P_0(dx)\\
	&=\IZ H\big(P(\sbt|X_0=x)\big|R(\sbt|X_0=x)\big)\, P_0(dx).
\end{align*}
Putting everything together, since $H(P_0|\mm)\le H(P|R)$ is finite,
\begin{align}\label{eq-19b}
\begin{split}
H(P &|R^{ P_0})=H(P|R)-H(P_0|\mm)\\
	&=\ud \big( H(P_T|\mm)-H(P_0|\mm)\big)+\ud E _{ P} \Iii \big(  |\bf ^{ P|R}|^2_\aa +|\bb ^{ P|R}|^2_\aa \big)(\Xb_t)/2\,dt.
\end{split}
\end{align}
From $ \left\{ \begin{array}{lcl}
\bcu&:=&(\bf -\bb)/2\\
\bos&:=&(\bf+\bb)/2
\end{array}\right.,$ we derive the parallelogram identity 
$$|\bf|^2_\aa /2+|\bb|^2_\aa /2=|\bcu|^2_\aa +|\bos|^2_\aa ,$$ leading to
\begin{align*}
H(P |R^{ P_0})
	= \mathcal{F}(P_T)- \mathcal{F}(P_0)+ E _{ P} \Iii \big(  |\bc { P|R}|^2_\aa +|\bo { P|R}|^2_\aa \big)(\Xb_t)/2\,dt.
\end{align*}
We conclude with Theorem \ref{res-08} and $H(P _{ [0,t]}|R ^{ P_0} _{ [0,t]})\le H(P|R ^{ P_0})< \infty,$ for all $0\le t\le T.$
\end{proof}

\subsection*{Heat flow}

In this subsection, the reference measure $R$ is defined with $T= \infty$, that is on $\OO=C([0, \infty),\ZZ).$

\begin{definition}
 The time marginal flow
$	
\mu_t:= R ^{ \mu_0}_t,  t\ge 0,
$	
of $R ^{  \mu_0}$
is called the \emph{heat flow} issued from $ \mu_0$.

\end{definition}

Next result is a direct consequence of Proposition \ref{res-12} which tells us that  the Fisher information is proportional to the rate of consumption of free energy along  the heat flow.
\begin{corollary}
If $H( \mu_0|\mm)< \infty,$ the heat flow $( \mu_t) _{ t\ge 0}$ satisfies
\begin{align*}
\mathcal{F}( \mu_t)- \mathcal{F}( \mu_0)
	=-2\int _{ [0,t]}  \mathcal{I}_\aa( \mu_s|\mm)\ ds,\quad \forall t\ge 0,
\end{align*}
where all these quantities are finite.
\end{corollary}

\begin{proof}
We have $0\le H( \mu_t|\mm)=H(R ^{ \mu_0}_t|R_t)\le H(R ^{ \mu_0}|R)=H( \mu_0|\mm)< \infty.$
Applying \eqref{eq-19b} with $P=R ^{ P_0= \mu_0}$  leads to 
\begin{align*}
0=H(R ^{ P_0} _{ [0,t]}|R ^{ P_0} _{ [0,t]})
	&=\mathcal{F}( \mu_t)- \mathcal{F}( \mu_0)
	+ E _{ R ^{ P_0}} \Iii \big(  |\bf ^{ R ^{ P_0}|R}|^2_\aa +|\bb ^{ R ^{ P_0}|R}|^2_\aa \big)(\Xb_t)/4\,dt\\
	&=\mathcal{F}( \mu_t)- \mathcal{F}( \mu_0)
	+\int _{ [0,t]\times \ZZ}  | \nabla \log \sqrt{ d \mu_s/d\mm}|^2_\aa \ dsd \mu_s,
\end{align*}
because $\bf ^{ R ^{ P_0}|R}=0$ implies that $\bb ^{ R ^{ P_0}|R}=2 \bo { R ^{ P_0}|R}=\nabla \log { d \mu_s/d\mm}.$
\end{proof}

\section{Random walks} \label{sec-rw}

In this section, the IbP formula  is used to obtain at Theorem \ref{res-21} a time reversal formula for a random walk on a graph under a finite entropy condition. This simple setting permits us to introduce Föllmer's guideline to derive time reversal formulas, with minimal technicalities.

\subsection*{Graph}

We consider continuous-time random walks on a countable graph $(\XX,\sim)$ where $\XX$ is the set of all vertices and the symmetric binary relation $x\sim y,$ $x,y\in\XX$ states that $\{x,y\}$ is a non-oriented edge of the graph. We assume without loss of generality that the graph is irreducible: $\XX$ is the unique class of communication, and that it contains no elementary loop: $x\sim x$ is forbidden.
We also assume that $(\XX,\sim)$ is a locally finite graph meaning that each vertex $x\in\XX$ admits finitely many neighbours. That is
\begin{equation}\label{eqRW-07}
n_x:= \# \left\{y\in\XX; y\sim x\right\} <\infty,\quad \forall x\in\XX.
\end{equation}
The countable set $\XX$ is equipped with its discrete topology.

\subsection*{Random walk}

A random walk on the graph $(\XX,\sim)$ is a time-continuous Markov measure $Q\in\MO$ which is specified by its initial distribution $Q_0\in\MX$ and  its forward generator acting on any real function  in the class 
\begin{align*}
\UU:= \left\{ u:\XX\to\RR; \#\supp(u)< \infty\right\} 
\end{align*}
 of all real functions  with a finite support via the formula
\begin{equation}\label{eqRW-01}
 \LLf^Q_t u(x)=\sy[u(y)-u(x)]\,\overrightarrow{j}(t,x;y), 
\quad  x\in\XX, t\in\ii,\quad u\in \mathcal{U},
\end{equation}
where for any adjacent neighbours $x\sim y,$ $\overrightarrow{j}(t,x;y)\ge 0$ is the average frequency of jumps from $x$ to  $y$ at time $t$. The jump kernel associated with this generator  is $$\sy \overrightarrow{j}(t,x;y)\delta_y\in\MX,\ x\in\XX, t\in\ii$$ where $\delta_y$ stands for the Dirac measure at $y$. 
For any pair of functions $u,v$ in $ \UU,$ the carré du champ  is
$$
 \Gaf^Q_t(u,v)(x)=\sy[u(y)-u(x)][v(y)-v(x)]\,\overrightarrow{j}(t,x;y).
 $$ 
Note  that the class of functions $ \mathcal{U}$ is  an algebra. See Remark \ref{rem-04} for the significance of this property.

\subsection*{A first time reversal formula}
We start by applying part (b) of Theorem \ref{res-02} as a first step of a more general result.

\begin{proposition}\label{res-27}
If for any $x\sim y$ the function $t\mapsto \overrightarrow{j}(t,x;y)$ is continuous, then
$\UU\subset \dom\LLf^Q\cap\dom\LLb^Q$,  
 the backward  generator is
\begin{equation*}\label{eqRW-01b}
 \LLb^Q_t u(x)=\sy[u(y)-u(x)]\,\overleftarrow{j}(t,x;y), 
\quad  x\in\XX, t\in\ii,\quad u\in \mathcal{U},
\end{equation*}
where for all $t\in\ii$  and all $ x,y\in\XX, x\neq y,$
\begin{align*}
\qq_t(x)\overleftarrow{j} (t,x;y)=\qq_t(y)\overrightarrow{j} (t,y;x).
\end{align*}
\end{proposition}

\begin{proof} 
Let us apply Theorem \ref{res-02}-(b).  
Under the assumptions \eqref{eqRW-07} and
\[\Iii \overrightarrow{j}(t,x;y)\,dt< \infty,\qquad \forall x,y: x\sim y,\]
the processes $\LLf^Qu(X)$ and $\Gaf^Qu(X)$ are in $L^1(\Qb),$ and
 \begin{align*}
 M ^{ Q,[u,v]}_t
 	= \sum _{0\le s\le t} [u(X_s)-u(X _{ s^-})][v(X_s)-v(X _{ s^-})]
		-\int _{ [0,t]} \Gaf^Q_s(u,v)(X_s)\,ds.
 \end{align*}
With our notation, this means that: $\UU^Q=\UU.$
\\
The class $\UU^Q=\UU$ determines the weak convergence of measures and our assumption about the continuity of $\overrightarrow{j}$ implies \eqref{eq-39b}. Denoting
$\overrightarrow{\qq j}(t,x;y):=\qq_t(x)\overrightarrow{j} (t,x;y)$, we see that for any $w\in\UU ^{ \Qb},$
\begin{multline*}
\IiX \Gaf^Q_t(u, w_t)[x]\,\qq_t(dx)dt
	=\Iii dt \sxy [u(y)-u(x)][w_t(y)-w_t(x)]\,\overrightarrow{\qq j}(t,x;y)\\
	=-\Iii dt \sxy w_t(x)[u(y)-u(x)]\,\big(\overrightarrow{\qq j}(t,x;y)
		+\overrightarrow{\qq j}(t,y;x)\big).
\end{multline*}
This proves that $\IiX \Gaf^Q(u, \sbt)\,d\qb$ is a finite measure, showing that the hypotheses of Theorem  \ref{res-02}-(b) are satisfied. Hence,    $u\in\dom\LLb^Q$,   $\LLb^Qu$ is integrable and for almost every $t$  the IbP formula \eqref{eq-24} holds, that is
\begin{multline*}
\sx v(x)\LLb^Q_tu(x)\, \qq_t(x)\\
= -\sxy  \{u(y)-u(x)\}v(x)\overrightarrow{\qq j}(t,x;y)
	-\sxy  \{u(y)-u(x)\}\{v(y)-v(x)\}\overrightarrow{\qq j}(t,x;y)\\
= -\sxy  \{u(y)-u(x)\}v(y)\overrightarrow{\qq j}(t,x;y)
= \sxy  \{u(y)-u(x)\}v(x)\overrightarrow{\qq j}(t,y;x),
\end{multline*}
for any $u,v\in\UU.$ 
\\
On the other hand, with \cite[Proposition 3.4]{Leo20}  we know that for almost all $t$ and for every $x$,
\begin{align}\label{eq-49}
\begin{split}
\LLb^Q_tu(x)
	&=\Limh h ^{ -1} E_Q \left[u(X _{ t-h})-u(X_t)\mid X_t=x\right]\\
	&= \sum _{ y\in\XX} \{u(y)-u(x)\} \Limh h ^{ -1}  Q(X _{ t-h}=y\mid X_t=x),
\end{split}
\end{align}
proving that the backward generator writes as  $\LLb^Q_tu(x) =\sum _{ y\in\XX}  \{u(y)-u(x)\}\,\overleftarrow{j} (t,x;y)$ for some function $\overleftarrow{j}.$ Plugging this into the expression $\sx v(x)\LLb^Q_tu(x)\, \qq_t(x),$ we arrive at 
\begin{align*}
\sxy \{u(y)-u(x)\}v(x) \qq_t(x)\overleftarrow{j} (t,x;y)
=\sxy  \{u(y)-u(x)\}v(x)\overrightarrow{\qq j}(t,y;x),
\end{align*}
and conclude remarking that the family of functions $(x,y)\mapsto  \{u(y)-u(x)\}v(x)$ when $u$ and $v$ describe $\UU$ is measure-determining off the diagonal of $\XXX.$
\end{proof}

\subsubsection*{Reversible random walk}

Saying that $Q\in\MO$ is reversible  means that there is a (possibly unbounded) positive measure $\mm\in\MX$ on $\XX$ such that, not only $Q$ is $\mm$-stationary i.e.\ 
$\qq_t=\mm,$ $  \forall 0\le t\le T,$ but also that $Q$ is invariant with respect to time reversal i.e.: for any subinterval $[r,t]\in\ii,$ 
$$(X _{(r+t-s)^-};r\le s\le t)\pf Q=(X_s;r\le s\le t)\pf Q.$$
This implies that the forward and backward transition mechanisms  do not depend on the time variable $t$ and are the same:  $\overrightarrow{j}=\overleftarrow{j}=:j$.
In view of Proposition \ref{res-27}, we obtain the detailed balance condition
\begin{equation}\label{eqRW-04}
\mm(x) j_x(y)=\mm(y) j_y(x),\quad \forall x,y\in\XX: x\sim y.
\end{equation}
Without loss of generality, we assume that $x\sim y \iff j_x(y),j_y(x)>0$ and that the graph is irreducible. It follows that $\mm(x)>0$ for all $x\in\XX.$ 
The general solution of \eqref{eqRW-04} is 
\begin{align*}
j_x(y)= s(x,y) \sqrt{\mm(y)/\mm(x)}
\end{align*}
where $s$ is any symmetric function such that  $x\sim y \iff s(x,y)>0.$


\subsubsection*{Counting  random walk}
If the waiting time at $x$ is  an exponential random variable $\mathcal{E}(n_x)$ and the jump occurs uniformly onto each neighbour, we obtain the jump kernel
\begin{equation}\label{eqRW-10}
J^o_x:= \sy \delta_y,\qquad x\in\XX,
\end{equation}
which
admits the counting measure 
\begin{equation}\label{eqRW-08}
\mm^o=\sx \delta_x\in\MX
\end{equation}
as a reversing measure. We denote $R\in\MO$ this reversible random walk with $R_0=\mm^o$ and call it the counting random walk. It will be the reference path measure for the rest of this section.

The remainder of this section is devoted to the proof of an extension of Proposition \ref{res-27}, stated at Theorem \ref{res-21}, where the hypothesis on the continuity of $\jf$  is removed and replaced by a finite entropy assumption.

\subsection*{Finite entropy assumption}

Let $P\in\PO$ be  a Markov probability measure such that $$H(P|R)< \infty,$$ with $R\in\MO$ the counting random walk. This finite entropy property implies (Girsanov's theory) that   there exists some measurable function 
$\jf^P:\iXX\to[0,\infty)$ which is defined $dtP_t(dx)J^o_x(dy)$-almost everywhere such that $P$ is the unique solution of the { martingale problem $ \textrm{MP}(P_0, \Jf^P)$}  associated to the initial marginal $P_0$ and the jump kernel  $\Jf^P=\jf^PJ^o,$ that is 
\begin{align*}
\Jf^P _{ t,x}=\sy \jf^P_{t,x}(y)\, \delta_y,\quad (t,x)\in\iX.
\end{align*}
Moreover, 
\begin{equation}\label{eqRW-12}
H(P| R)=H(P_0|\mm^o)+\IiXX \hh\Big(\jf^P_{t,x}(y)\Big)\,dt\pp_t(dx)J^o_x(dy)<\infty,
\end{equation}
where 
\begin{equation*}
\hh(a):=\left\{\begin{array}{ll}
a\log a-a+1, & \textrm{if } a>0,\\
1,& \textrm{if }a=0,\\
\infty,& \textrm{if }a<0.
\end{array}\right.
\end{equation*}

\begin{lemma}\label{res-17}
Let $u$ be any function in $\UU$ and $P\in\PO$ satisfy $H(P|R)< \infty.$

\begin{enumerate}[(a)]
\item
The function $u$ stands in $\dom\LLf^P$, $\LLf^Pu(\Xb)\in L\log L (\Pb)$ and 
$$\LLf^P _tu(x)=\sy [u(y)-u(x)]\ \jf^P _{ t,x}(y).$$
\item
There exists some measurable function 
$\jb^P:\iXX\to[0,\infty)$ which is defined $dt\pp_t(dx)J^o_x(dy)$-almost everywhere such that $u$ stands in $\dom\LLb^P$, with $\LLb^Pu(\Xb)\in L\log L (\Pb)$ and 
\begin{align}\label{eq-27}
\LLb^P _tu(x)=\sy [u(y)-u(x)]\ \jb^P _{ t,x}(y).
\end{align}
Moreover
\begin{equation*}
H(P| R)=H(P_T|\mm^o)+\IiXX \hh\Big(\jb^P_{t,x}(y)\Big)\,dt\pp_t(dx)J^o_x(dy)<\infty.
\end{equation*}
\end{enumerate} 
\end{lemma}

\begin{proof}
\boulette{(a)}
Let us denote the right hand side of the desired identity by: $A(t,x):=\sy [u(y)-u(x)]\ \jf^P _{ t,x}(y).$
 With \eqref{eqRW-12}, our assumption \eqref{eqRW-07}, the finiteness of the support of $u,$ and
\begin{align*}
\begin{split}
| A(t,x)|
	\le 2 \sup |u| \1 _{ \left\{x\in\supp(u)\right\} }  \sy     \jf^P _{ t,x}(y),
\end{split}
\end{align*}
we see  that $A(\Xb)$ is in $L\log L (\Pb)$. This implies that $\Iii|A(\Xb_t)|\,dt$ is finite $P\ae$, so that $A(t,x)=\LLf^P u_t(x),$ $\pb(dtdx)\ae$ and $\LLf^Pu(\Xb)\in L\log L (\Pb).$

\Boulette{(b)}
Time reversal being a bijective mapping: $H(P^*|R^*)=H(P|R),$ see Proposition \ref{res-16}. Since $R$ is chosen to be reversible, we also have $R^*=R,$ leading to: $$H(P^*|R)=H(P|R)< \infty.$$ Hence we are allowed to apply (a) which  tells us that $u\in\dom\LLf ^{ P^*},$ $\LLf ^{ P^*}u(\Xb)\in L\log L (\overline{P^*})$, and there is some measurable function $\jf ^{ P^*}$ such that 
$\LLf ^{ P^*} _tu(x)=\sy [u(y)-u(x)]\ \jf ^{ P^*} _{ t,x}(y).$ We conclude taking $\jb^P _{ t,x}(y):=\jf ^{ P^*} _{ T-t,x}(y).$
\end{proof}

\subsection*{Time reversal formula}
 The main theorem of this section is

\begin{theorem}\label{res-21}
Let $P\in\PO$ be a Markov random walk with forward generator
\begin{align*}
\LLf^P _tu(x)=\sy [u(y)-u(x)]\ \jf^P _{ t,x}(y),
\qquad x\in\XX,\ u\in\UU,
\end{align*}
where the forward intensity of jump $\jf^P$ is measurable.
If $H(P|R)< \infty,$ i.e.\ 
\begin{align*}
H(\pp_0|\mm^o)+\Iii dt\sxy \hh\Big(\jf^P_{t,x}(y)\Big)\, \pp_t(x)<\infty,
\end{align*}
then, $\UU\subset\dom\LLb^P$ and 
\begin{align*}
\LLb^P _tu(x)=\sy [u(y)-u(x)]\ \jb^P _{ t,x}(y),
\qquad x\in\XX,\ u\in\UU,
\end{align*} 
where  the backward intensity of jump $\jb^P$ verifies for almost all $0\le t\le T,$
\begin{align}\label{eq-51}
\pp_t(x)\jf^P_{t,x}(y)=\pp_t(y)\jb^P _{ t,y}(x),\quad 
\forall x, y\sim x \in\XX.
\end{align}
\end{theorem}

\begin{proof}
With \eqref{eqRW-12} we see that $H(P|R)< \infty,$ i.e.\ the hypothesis of  Lemma \ref{res-17} is satisfied. This lemma tells us that   the assumptions of Theorem \ref{res-02}-(a) are satisfied. Therefore, for almost all $t\in\ii,$  and any $u,v\in\UU,$ the IbP formula \eqref{eq-24} holds, i.e.
\begin{eqnarray*}
0&=& \IXX [u(y)-u(x)]v(x)(\jf^P+\jb^P)_{t,x}(y)\,\pp_t(dx)J^o_x(dy)\\
&&\quad +\IXX [u(y)-u(x)][v(y)-v(x)]\jf^P_{t,x}(y)\,\pp_t(dx)J^o_x(dy)\\
&=& \IXX [u(y)-u(x)][v(x)\jb^P_{t,x}(y)+v(y)\jf^P_{t,x}(y)] \pp_t(x)\,\mm^o(dx)J^o_x(dy).
\end{eqnarray*}
As the counting random walk is $\mm^o$-reversible, by \eqref{eqRW-04} $\mm^o(dx)J^o_x(dy)$ is a symmetric measure on $\XXX$ (obvious by direct inspection). It follows that
\begin{equation*}
0=\IXX [u(y)-u(x)]v(y)[ \pp_t(x)\jf^P_{t,x}(y)- \pp_t(y)\jb^P _{ t,y}(x)]\,\mm^o(dx)J^o_x(dy),
\end{equation*}
from which the result follows.
\end{proof}

\appendix

\section{Stochastic derivatives and extended generators} \label{sec-leo}

After recalling the definitions of Markov measures, extended generators and stochastic derivatives, we state a couple of technical results obtained in
 \cite{Leo20}.
 
 \begin{definition}[Markov measure] \label{def-01}
 A path measure $Q$ such that  $Q_t$ is  $ \sigma$-finite for all $t$ is called a conditionable path measure.
  A path measure $Q\in\MO$ is said to be Markov if it is conditionable and for any $0\le t\le T,$ $Q(X _{ [t,T]}\in\sbt\mid X _{ [0,t]})=Q(X _{ [t,T]}\in\sbt\mid X _t).$ 
 \end{definition}

  The reason for requiring $Q$ to be conditionable is that  it allows for defining  the conditional expectations $E_Q(\sbt\mid X_ \mathcal{T})$ for any $ \mathcal{T}\subset \ii$ even in the case where $Q$ is an unbounded measure, see \cite[Def.\,1.10]{Leo12b}.

 The notion of extended generator was introduced by H. Kunita \cite{K69} and extensively used by P.\,A.\,Meyer and his collaborators, see \cite{DM4}. Here is a variant of this definition.

\begin{definition}[Extended forward generator of a Markov measure]\label{def-02}
Let $Q$ be a Markov measure. A measurable function $u$ on $\iX$ is said to be in
the domain of the extended forward generator of $Q$ if there exists a
measurable function $v$ on $\iX$  such that
 $\Iii|v(\Xb_t)|\,dt<\infty,$ $Q\ae$ and the process $$M^u_t:=u(\Xb_t)-u(\Xb_0)-\II0t v(\Xb_s)\,ds,\quad 0\le t\le T,$$ is a local $Q$-martingale.
We denote
$$
   \LLf^Q u(t,x):=v(t,x)
$$
and call $\LLf^Q $ the extended forward generator of $Q.$ Its domain  is denoted by $\dom
\LLf^Q .$
\end{definition}

\begin{remarks}\label{rem-02}\ \begin{enumerate}[(a)]
 \item
 In the case where $Q$ is the law of a Markov process associated with some semigroup with generator $ \mathcal{G}$ and  $u:\iX\to\RR$ is a $t$-differentiable function  such that for each $t,$ $u(t,\sbt)$ belongs to the domain of $ \mathcal{G},$ then $u$ belongs to $\dom\LLf^Q$ and
 \begin{equation*}
 \LLf^Q u=(\partial_t+ \mathcal{G})u.
 \end{equation*}

\item
The notation $v=\mathcal{L}u$ almost rightly suggests that $v$ is
a function of $u.$ Indeed, when $u$ is in $\dom\LLf^Q ,$
the Doob-Meyer decomposition of the special semimartingale
$u(\Xb_t)$  into its predictable bounded variation part $\int
v_s\,ds$ and its local martingale part is unique. But one can
modify $v=\LLf^Q u$ on a small (zero-potential) set without
breaking the martingale property. As a consequence, $u\mapsto
\LLf^Q u$ is a multivalued operator and $u\mapsto
\LLf^Q u$ is an almost linear operation.
 \end{enumerate}

\end{remarks}

Extended generators are connected with   martingale problems which were introduced by Stroock and Varadhan \cite{SV79}.

\begin{definition}[Martingale problem]
Let $\mathcal{C}$ be a class of measurable real functions $u$ on $\iX$
and for each $u\in\mathcal{C},$ let $\mathcal{L}u:\iX\to\RR$  be a
measurable function such that $\Iii |\mathcal{L}u(t,\omega_t)|\,dt<\infty$ for all $\omega\in\OO.$ Take also a positive $ \sigma$-finite measure
$\mu_0\in\MX.$ One says that $Q\in\MO$ is a solution to the
martingale problem $\MP(\mathcal{L},\mathcal{C};\mu_0)$ if
$Q_0=\mu_0\in\MX$ and for all $u\in\mathcal{C},$ the process
$
    u(\Xb_t)-u(\Xb_0)-\II0t \mathcal{L}_su[X_s]\,ds
$
is  a local $Q$-martingale.
\end{definition}
 
 Proposition \ref{res-10} below states that the extended generator can be computed by means of a stochastic derivative. 
Nelson's definition \cite{Nel67} of the stochastic derivative is the following.

\begin{definition}[Stochastic forward derivative of a Markov measure]\label{def-03}
Let $Q$ be a Markov measure and $u$ be a measurable real function on $\iX$ such that $E_Q|u(\Xb_s)|<\infty$ for all $0\le s\le T.$ 
\begin{enumerate}
\item
We say
that $u$ admits a stochastic forward derivative  under $Q$ at time
$t\in[0,T)$ if   the following limit
\begin{align*}
 \Lf^Qu(t,x):=\Limh E_Q\left(\frac1h [u(\Xb_{t+h})-u(t,x)]
    \mid  X_t=x  \right)
\end{align*}
exists in $L^1(\ZZ,Q_t).$
\\
In this case,  $\Lf^Qu(t,\sbt)$ is called  the stochastic forward derivative of $u$ at time $t$.
\item
If $u$ admits a stochastic forward derivative for  almost all
$t,$ we say that $u$ belongs to the domain $\dom \Lf^Q $ of
the stochastic forward derivative $\Lf^Q$ of $Q.$
\item
If $u$ does not depend on the time variable $t,$ we
denote
$
    \Lf^Q_tu(x):=\Lf^Qu(t,x).
$
\end{enumerate}
\end{definition}

\subsubsection*{Reversing time}

If $Q\in\MO$    is Markov,  so is its time reversal $Q^*$, and one can consider the extended generators and stochastic derivatives of both $Q$ and $Q^*.$ More generally, we introduce the following notions.

As a notation, the $ \sigma$-field generated by $ X _{ [t^-,T]}$ is $ \sigma(X _{ [t^-,T]}):= \cap _{ h>0} \sigma(X _{ [t-h,T]})= \sigma(X _{ t^-})\vee \sigma(X _{ [t,T]}).$

\begin{definition}[Extended backward generator]\label{def-02b}
Let $Q$ be a conditionable path measure. A process $u$ adapted to the predictable backward filtration $( \sigma(X _{ [t^-,T]}); 0\le t\le T)$ is said to be in
the domain of the extended backward generator of $Q$ if there exists a
process $v$ also  adapted to the predictable backward filtration such that
 $\Iii|v(t,X _{ [t^-,T]})|\,dt<\infty,$ $Q\ae$ and the process 
 $$
u(t, X _{ [t^-,T]})-u(T,X_T)-\II tT v(s, X _{ [s^-,T]})\,ds,\quad 0\le t\le T,
 $$ 
 is a local backward $Q$-martingale.
We denote
$$
   \LLb^Q_t u:=v_t
$$
and call $\LLb^Q $ the extended backward generator of $Q.$ Its domain
  is denoted by $\dom
\LLb^Q .$
\end{definition}

\begin{definition}[Stochastic backward derivative]\label{def-03b}
Let $Q$ be a conditionable path measure and a measurable function $u$ on $\iX$  such that $E_Q|u(s, X _s)|<\infty$ for all $0\le s\le T.$ 
\begin{enumerate}
\item
We say
that $u$ admits a stochastic backward derivative  under $Q$ at time
$t\in(0,T]$ if   the following limit
\begin{align*}
 \Lb^Qu(t,X _{ [t^-,T]}):=\Limh E_Q\left(\frac1h [u(\Xb_{t-h})-u(\Xb_t)]
    \mid  X _{ [t^-,T]}  \right)
\end{align*}
if this limit exists in $L^1(Q).$  
\\
In this case,  $\Lb^Qu(t,\sbt)$ is called  the stochastic backward derivative of $u$ at time $t$.
\item
If $u$ admits a stochastic backward derivative for  almost all
$t,$ we say that $u$ belongs to the domain $\dom \Lb^Q $ of
the stochastic backward derivative $\Lb^Q$ of $Q.$
\end{enumerate}
\end{definition}

\subsection*{Convergence results}

A useful technical result for our purpose  is the following  convolution result.

\begin{lemma}\label{res-09}
For all $h>0,$ let
$k^h$ be a measurable nonnegative  convolution kernel such that $\supp
k^h\subset[-h,h]$ and $\int_{\mathbb{R}}k^h(s)\,ds=1.$
Let $Q$ be a $ \sigma$-finite positive measure on $\OO$ and  $v$ be a process in $L^p(\Qb)$ with    $1\le p< \infty$. 
\\
Define for all $h>0, t\in\ii$ and $\omega\in\OO$,\   $k^h *
v(t, \omega):=\int_{\ii}k^h(t-s)v_s( \omega)\,ds$.\\
Then, $k^h *v$ is in $L^p(\Qb)$ and
    $
    \Limh k^h *v=v\  \textrm{in }L^p(\Qb).
    $
\end{lemma}

We see that $k^h (s)\, ds$ is a probability measure on $\RR$ which
converges narrowly to the Dirac measure $\delta_0$ as $h$ tends down to zero. 
 We shall  invoke this lemma with $p=1$ or $2.$

\begin{corollary}\label{res-14}
Assume that in addition to the hypotheses of Lemma \ref{res-09}, for any $0\le t\le T,$ the random variable $v(t,\sbt)$ is $  \mathcal{A}_t$-measurable where $ \mathcal{A}_t$ is some sub-$ \sigma$-field.  Then,  the process $v^h$  defined by $v^h_t:= E_Q[k^h*v(t)\mid \mathcal{A}_t],$  is in  $L^p(\Qb)$ and $\Limh v^h=v$ in $L^p(\Qb).$
\end{corollary}

\begin{proof}
By Jensen's inequality
\begin{align*}
\|v^h-v\| _{ p,\Qb}^p
	&=\int _{ \XXb} |E_Q[k^h*v(t) \mid \mathcal{A}_t]-v(t)|^p\,d\Qb
	=\int _{ \XXb} |E_Q[k^h*v(t)-v(t) \mid \mathcal{A}_t]|^p\,d\Qb\\
	&\le \int _{ \XXb} E_Q[|k^h*v(t)-v(t)|^p \mid \mathcal{A}_t]\,d\Qb
	=\int _{ \XXb} E_Q|k^h*v(t)-v(t)|^p  \,d\Qb\\
	&=\|k^h*v-v\| _{p, \Qb}^p \underset{h\to 0^+}\longrightarrow 0,
\end{align*}
where the vanishing limit is the content of Lemma \ref{res-09}.
\end{proof}

Next proposition states that extended generators and stochastic derivatives are essentially the same.

\begin{proposition}\label{res-10} Let $Q$ be a conditionable measure.

\begin{enumerate}[(a)]
\item
 If  $u$ is in $\dom\LLf^Q$ and satisfies
$E_Q\Iii
 \big|\LLf^Q u(t,X _{ [0,t]})\big|^p\,dt<\infty$ for some  $p \ge 1,$ then
\begin{equation*}
    \Limh E_Q\II0{T-h} \left|
    \frac 1hE_Q\Big[u(\Xb_{t+h})-u(\Xb_t)\mid X _{ [0,t]} \Big]-\LLf^Q u(t,X _{ [0,t]})\right|^p\,dt=0.
\end{equation*}
In particular, this implies that
$u\in\dom \Lf^Q,$ and
the limit 
$$
\LLf^Q u(\sbt, X _{ [0, \sbt]})=\Lf^Q u(\sbt, X _{ [0, \sbt]}):=\Limh \frac 1hE_Q\Big[u(\Xb _{ \sbt+h})-u(\Xb_{\sbt})\mid X _{ [0, \sbt]} \Big]
$$ 
takes place in $L ^p(\Qb).$

\item
If  $u$ is in $\dom\LLb^Q$ is such that
$E_Q\Iii
 \big|\LLb^Q u(t, X _{ [t,T]})\big|^p\,dt<\infty$ for some  $p \ge 1,$ then
\begin{equation*}
    \Limh E_Q\II h{T} \left|
    \frac 1hE_Q\Big[u(\Xb_{t-h})-u(\Xb_t)\mid X _{ [t^-,T]} \Big]-\LLb^Q u(t,X _{ [t^-,T]})\right|^p\,dt=0.
\end{equation*}
In particular, this implies that
$u\in\dom \Lb^Q,$ and
the limit 
$$\LLb^Q u(\sbt, X _{ [\sbt,T]})=\Lb^Q u(\sbt, X _{ [\sbt,T]}):=\Limh \frac 1hE_Q\Big[u(\Xb _{\sbt -h})-u(\Xb _{ \sbt})\mid X _{ [\sbt,T]} \Big]$$ takes place in $L ^p(\Qb).$
\end{enumerate}
\end{proposition}

\begin{proposition}\label{resh-03b}\ 
\begin{enumerate}[(a)]
\item

Let  $u$ be  a  measurable real function on $\XXb$, and  $v$ be a  forward-adapted process
such that $u(\Xb)$ and $v$ are $\Qb$-integrable,  $t\mapsto u(\Xb_t)$ is right continuous (for instance $u$ might be continuous) and
 \begin{equation}\label{eqh-03}
    \Limh E_Q\II0{T-h}\left|\frac1h E_Q[u(\Xb_{t+h})-u(\Xb_{t})\mid X _{ [0,t]}]-v_t\right|\,dt =0.
\end{equation}
Then, $u$ belongs to $\dom\LLf^Q$ and $\dom L^Q$,  and
    $
      \LLf^Q u=L^Q u=v,\ \Qb\ae
    $
    
\item
Let $u$ be a  measurable real function on $\XXb$ and $v$ a  backward-predictable process,
such that $u(\Xb), v$ are $\Qb$-integrable,  $t\mapsto u(\Xb^*_t)$ is right continuous (for instance $u$ might be continuous) and
 \begin{equation*}
    \Limh E_Q\II h{T}\left|\frac1h E_Q[u(\Xb_{t-h})-u(\Xb_{t})\mid X _{ [t^-,T]}]-v(t, X _{ [t^-,T]})\right|\,dt =0.
\end{equation*}
Then, $u$ belongs to $\dom\LLb^Q$ and $\dom \Lb^Q$,  and
    $
      \LLb^Q u=\Lb^Q u=v,\ \Qb\ae
    $

\end{enumerate}
\end{proposition}

\section{Relative entropy with respect to an unbounded measure} \label{sec-um}

Let $r$ be some $\sigma$-finite positive measure on some measurable space $Y$. The relative entropy of the probability measure $p$ with respect to $r$ is loosely defined by
$$
H(p|r):=\int_Y \log(dp/dr)\, dp\in (-\infty,\infty],\qquad p\in \mathrm{P}(Y)
$$
if $p\ll r$ and $H(p|r)=\infty$ otherwise. More precisely, when $r$ is a probability measure,  we have $H(p|r)=\int_Y h(dp/dr)\,dr\in[0,\infty]$ with $h(a)=a\log a-a+1\ge 0$ for all $a\ge0,$ (take $h(0)=1).$ Hence  this definition is meaningful and it follows from the strict convexity of $h$ that $H(\sbt|r)$ is also strictly convex.
\\
If $r$ is unbounded, one must restrict the definition of $H(\sbt|r)$ to some subset of $\mathrm{P}(Y)$ as follows. As $r$ is assumed to be $\sigma$-finite, there exist  measurable functions $W:Y\to [0,\infty)$ such that
\begin{equation}\label{eqRW-05}
z_W:=\int_Y e ^{-W}\, dr<\infty.
\end{equation}
Define the probability measure $r_W:= z_W ^{-1}e ^{-W}\,r$ so that $\log(dp/dr)=\log(dp/dr_W)-W-\log z_W.$ It follows that for any $p\in \mathrm{P}(Y)$ satisfying $\int_Y W\, dp<\infty,$ the formula 
\begin{equation*}
H(p|r):=H(p|r_W)-\int_Y W\,dp-\log z_W\in (-\infty,\infty]
\end{equation*}
is a meaningful definition of the relative entropy which is coherent in the following sense. If $\int_Y W'\,dp<\infty$ for another measurable function $W':Y\to[0,\infty)$ such that $z_{W'}<\infty,$ then $H(p|r_W)-\int_Y W\,dp-\log z_W=H(p|r _{W'})-\int_Y W'\,dp-\log z_{W'}\in (-\infty,\infty]$.
\\
Therefore, $H(p|r)$ is well-defined for any $p\in \mathrm{P}(Y)$ such that $\int_Y W\,dp<\infty$ for some measurable non-negative function $W$ verifying \eqref{eqRW-05}.

It is well known that the relative entropy with respect to a \emph{probability} measure $r$ is invariant with respect to the push-forward by an injective  mapping. This is still true if $r$ is unbounded.

\begin{proposition}\label{res-16}
Let $r$  and $H(\sbt|r)$ be as above, and let $f:Y\to Z$ be a measurable mapping. For any $p\in \mathrm{P}(Y)$ satisfying $\int_Y W\,dp< \infty,$ we have: 
$H(f\pf p| f\pf r)\le H(p|r).$
\\
If in addition $f$ is injective, then $H(f\pf p| f\pf r)= H(p|r).$
\end{proposition}

\begin{proof}
It is a direct consequence of the variational formula
\begin{align*}
H(p|r)
	= \sup _{ u\in B_W(Y)} \left\{ \int_Y u\,dp-\int_Y e ^{ u-1}\,dr  \right\},
\end{align*}
where $B_W(Y):= \left\{u:Y\to\RR, \sup_Y |u|/(1+W) < \infty\right\}.$ Indeed
\begin{align*}
H(f\pf p|f\pf r)
	&= \sup _{ v \in B _{ W\circ f ^{ -1}}(f(Y))} \left\{ \int_{f(Y)} v\,d(f\pf p)-\int_{f(Y)} e ^{ v-1}\,d(f\pf r)  \right\}\\
	& = \sup _{  v \in B _{ W\circ f ^{ -1}}({f(Y)})} \left\{ \int_Y v\circ f\,dp-\int_Y e ^{ v\circ f-1}\,dr  \right\}\\
	&\le  \sup _{ u\in B_W(Y)} \left\{ \int_Y u\,dp-\int_Y e ^{ u-1}\,dr  \right\}
	 =H(p|r),
\end{align*}
because $ \left\{ v\circ f; v\in B _{ W\circ f ^{ -1}}({f(Y)}) \right\} \subset B_W(Y).$
\\
If $f$ is injective, 
 $u=v\circ f$ describes $B_W(Y)$ when $v$ describes $B _{ W\circ f ^{ -1}}({f(Y)})$, leading to an equality. 
\end{proof}


\end{document}